%% file: block-gmres-stag.tex
\documentclass{siamart0516}
\pdfoutput=1

\usepackage{amsfonts}
\usepackage{graphicx}
\usepackage{epsfig}
\usepackage{euscript,mathtools,amssymb}
\usepackage{enumerate}
\ifpdf
  \DeclareGraphicsExtensions{.eps,.pdf,.png,.jpg}
\else
  \DeclareGraphicsExtensions{.eps}
\fi
\usepackage{xcolor}
\definecolor{darkgreen}{rgb}{0.0,0.6,0.0}

\newcommand{\TheTitle}{Stagnation of block GMRES and its relationship to block FOM} 
\newcommand{\TheAuthors}{Kirk M. Soodhalter}
\headers{\TheTitle}{\TheAuthors}

\title{{\TheTitle}}

\author{
  Kirk M. Soodhalter\thanks{Johann Radon Institute for Computational and Applied Mathematics, Linz, Austria
    (\email{kirk.soodhalter@ricam.oeaw.ac.at}, \url{http://math.soodhalter.com}).}
}

\usepackage{amsopn}
\DeclareMathOperator{\diag}{diag}

\input{macros}

\begin{document}

\maketitle

\begin{abstract}
We analyze the the convergence behavior of block GMRES and 
characterize the phenomenon of stagnation which is then
related to the behavior of the block FOM method.  We 
generalize the block FOM method to generate well-defined approximations in the
case that block FOM would normally break down, and these generalized solutions
are used in our analysis.  This behavior is also related to the principal angles between the
column-space of the previous block GMRES residual  and the current minimum residual constraint space.
 At iteration $j$, it is shown that the proper generalization of GMRES stagnation to the
block setting relates to the columnspace of the $j$th block Arnoldi vector.
Our analysis covers both the cases of normal iterations as well as block Arnoldi breakdown wherein 
dependent basis vectors are replaced with random ones.
Numerical examples are given to illustrate what
we have proven, including a small application problem
to demonstrate the validity of the analysis in a less pathological case.
\end{abstract}

\begin{keywords}
Block Krylov subspace methods, GMRES, FOM, Stagnation
\end{keywords}

\begin{AMS}
65F10, 65F50, 65F08
\end{AMS}

\section{Introduction}\label{section.intro}
The Generalized Minimum Residual Method (GMRES) \cite{Saad.GMRES.1986} and 
the Full Orthogonalization Method (FOM) \cite{Saad1980} are two Krylov 
subspace methods for solving linear systems with non-Hermitian coefficient matrices and
one right hand side, i.e., 
\be\label{eqn.Axb}
	\vek A\vek x = \vek b\mwith \vek A\in\C^{n\times n}\mand \vek b\in\C^{n}.
\ee
The convergence behavior of these two methods is closely related, and this relationship
was characterized by Brown \cite{B.1991}, and other related results can be found in 
\cite{CG.1996,C.1994,W.1995}, and a related detailed geometric analysis of projection methods was
presented in \cite{Eiermann2001}.  A nice description can also be found
in \cite[Section 6.5.5]{Saad.Iter.Meth.Sparse.2003}.  Krylov subspace methods 
have been generalized to treat the situation in which we have multiple right-hand sides,
i.e., we are solving 
\be\label{eqn.AXB}
	\vek A\vek X = \vek B\mwith\vek B\in\C^{n\times L}.
\ee
In particular, block GMRES and block FOM \cite[Section 6.12]{Saad.Iter.Meth.Sparse.2003}
 have been proposed for solving 
\eqref{eqn.AXB}; however, to our knowledge, a similar full analysis of block GMRES, the connection between
stagnation and block FOM convergence, and accompanying geometric considerations have yet to be
described in the literature.  Therefore, in this work
we analyze the stagnation behavior of block GMRES and characterize its relationship
to the behavior of the block FOM method.  Similar analytic tools as in
in \cite{B.1991} and \cite{Eiermann2001} are used, but the behavior of block methods is a bit more complicated to describe.
The key result is the proper generalization of GMRES stagnation to the block setting.  The analog of stagnation
for block GMRES is not simply stagnation of some columns of the iterate.  Rather, at iteration $j$ it is associated to the dimension of the
intersection between the column space of the $j$th block Arnoldi vector and the $j$th block GMRES correction.  Stagnation of
some columns of the iterate is shown to be a special case thereof.  This then allows analogs of many of the
results on stagnation of GMRES and the relationship between GMRES and FOM to be proven in the block setting.
As block methods can suffer from partial or full stagnation of the iteration and breakdowns due to linear
dependence of the block residual, additional analysis 
is needed to fully characterize the stagnation in these settings. Here we consider the case that dependent basis vectors
are replaced with random ones (as in \cite{B.2000,CK.2010,Parks.Soodhalter.Szyld.16,S.2014}). One could similarly consider
the case that dependent vectors are removed and the block size reduced; see, e.g., 
\cite{A-Giraud-YanFei.GM-Inexact.2014, Yan-Fei.GCH.2016-inprep,RS.2006, OLeary1980}.

The rest of this paper proceeds as follows.  In the next section, we review 
Krylov subspace methods, focusing in particular on block GMRES and block FOM.
We also review existing analysis relating GMRES- and FOM-like methods.  The type of
relationship illuminated in \cite{B.1991} has been extended to many other pairs of 
methods.  In Section \ref{section.main-results}, we present our main results which 
characterize the relationship between block GMRES and block FOM.  In Section 
\ref{section.numerical-examples}, we construct numerical examples which demonstrate
what has been revealed by our analysis.  We offer some discussion and conclusions in
Section \ref{section.conclusions}.

In this paper, we adopt the convention that $\vek I$ is the identity matrix, where context determines the appropriate dimension.
When needed, we specify the dimension $\vek I_{J}\in\R^{J\times J}$.  Similarly, $\vek 0$ denotes the matrix of zeros, with
dimension determined by context.  We denote $\vek 0_{J}\in\R^{J \times J}$ to be a square matrix of zeros and 
$\vek 0_{J_{1}\times J_{2}}\in\R^{J_{1}\times J_{2}}$ with $J_{1}\neq J_{2}$ to be a rectangular matrix of zeros.

\section{Background}\label{section.prelim}
In this section, we review the basics about Krylov subspace methods and focus on the
the block version, designed to solve, e.g., \eqref{eqn.AXB}.  
We describe everything in terms of block Krylov subspace methods, and discuss the simplifications in the 
case that the block size $L=1$.
We then review existing
results relating the iterates of pairs of methods 
(many times derived from Galerkin and minimum
residual projections, respectively), 
e.g., FOM and GMRES \cite{B.1991} and
BiCG and QMR \cite{freund.QMR.1991} as well as subsequent works which expand upon
and offer additional perspective on these pair-wise relationships, e.g.,  
\cite{CG.1996,C.1994,GR.2001,RS.2002,W.1995}.
\subsection{Single-vector and block Krylov subspaces}\label{section.Krylov-review}

In the case that we are solving the system \eqref{eqn.AXB} with multiple
right-hand sides (a block right-hand side), block Krylov subspace methods are 
an effective family of methods for generating high quality approximate solutions 
to \eqref{eqn.AXB} at relatively low cost.
Let $\vek X_{0}$ be an initial approximate solution to \eqref{eqn.AXB} with
block initial residual $\vek F_{0}~=~\vek B - \vek A\vek X_{0}$. We can define the $j$th block Krylov subspace
as 
\be\label{eqn.block-Krylov-def}
	\K_{j}(\vek A,\vek F_{0}) = \Span\curl{\vek F_{0}, \vek A\vek F_{0}, \vek A^{2}\vek F_{0}, \ldots, \vek A^{j-1}\vek F_{0}}
\ee
where the span of a collection of block vectors is understood to be the span of all their 
columns.  When $L=1$ ($\vek B,\,\vek X_{0}\in\C^{n}$), this definition reduces to the single-vector
Krylov subspace, denoted $\CK_{j}(\vek A,\vek F_{0})$. 
In the case $L>1$, is straightforward to show that 
\be\nn
	\K_{j}(\vek A,\vek F_{0}) = \CK_{j}(\vek A, \vek F_{0}(:,1)) + \CK_{j}(\vek A, \vek F_{0}(:,2)) + \cdots + \CK_{j}(\vek A, \vek F_{0}(:,L))
\ee
where we use the MATLAB style indexing notation $\vek F(:,i)$ to 
denote the $i$th column of a matrix $\vek F\in\C^{I\times J}$ such that $J\geq i$; see, e.g., \cite{GS.2009}.

Let $\vek W_{j} = \bbmat \vek V_{1} & \vek V_{2} & \ldots \vek V_{j} \ebmat\in\C^{n\times jL}$ be the matrix with orthonormal columns spanning 
$\K_{j}(\vek A, \vek F_{0})$ with $\vek V_{i}\in\C^{n\times L}$ 
having orthonormal columns,
and $\vek V_{i}^{\ast}\vek V_{j}=\vek 0$ for $i\neq j$.  These orthonormal blocks can 
be generated one block at a time by an iterative orthogonalization process called 
the block Arnoldi process, which is a natural generalization of the
Arnoldi process for the single-vector case.  We have the block Arnoldi relation
\be\label{eqn.block-arnoldi-relation}
	\vek A\vek W_{j} = \vek W_{j+1}\overline{\vek H}_{j}^{(B)}
\ee 
where $\overline{\vek H}_{j}^{(B)}=\prn{\vek H_{i,j}}\in\C^{(j+1)L\times jL}$ is block 
upper Hessenberg with $\vek H_{ij}\in\C^{L\times L}$ and $\vek H_{j+1,j}$ upper 
triangular.  

We can derive block FOM and block GMRES
methods through Galerkin and minimization constraints. We have for each column of the $j$th
block residual the constraints
\bea
	\vek R_{j}(:,i)&\perp &\K_{j}(\vek A, \vek F_{0})\mor\label{eqn.block-galerkin-cond}\\ 
	\vek R_{j}(:,i)&\perp &\vek A\K_{j}(\vek A, \vek F_{0})\label{eqn.block-min-res-cond}.
\eea
which lead to the block FOM and block GMRES methods, respectively.  For both methods,
approximations can be computed for all columns simultaneously.
Let $\vek X_{j}^{(F)}$ and $\vek X_{j}^{(G)}$ denote the $j$th block FOM and block GMRES 
approximation solutions for \eqref{eqn.AXB}.  Furthermore, let 
$\vek E_{L}^{[I]}\in\R^{I\times L}$ have as columns the first $L$ columns of the
$I\times I$ identity matrix, and let $\vek F_{0}=\vek V_{1}\vek S_{0}$ be the reduced
QR-factorization with $\vek S_{0}\in\C^{L\times L}$ upper-triangular.
Using \eqref{eqn.block-arnoldi-relation}, block FOM can be derived
from \eqref{eqn.block-galerkin-cond} which leads to the formulation
\be\label{eqn.block-FOM-subprob}
	\vek X_{j}^{(F)} = \vek X_{0} + \vek T_{j}^{(F)}\mwhere \vek T_{j}^{(F)} = \vek W_{j}\vek Y_{j}^{(F)}\mand \vek H_{j}^{(B)}\vek Y_{j}^{(F)} = \vek E_{L}^{[jL]}\vek S_{0},
\ee
where $\vek H_{j}^{(B)}\in\C^{jL\times jL}$ is defined as the matrix containing the first $jL$ 
rows of $\overline{\vek H}_{j}^{(B)}$.
Similarly for block GMRES, we can use \eqref{eqn.block-arnoldi-relation}, combined with 
\eqref{eqn.block-min-res-cond} to yield a formulation
\bea
	\vek X_{j}^{(G)} = \vek X_{0} + \vek T_{j}^{(G)} &\mwhere & \vek T_{j}^{(G)} = \vek W_{j}\vek Y_{j}^{(G)}\nn \\ &\mand & \vek Y_{j}^{(G)}=\argmin{\vek Y\in\C^{jL\times L}}\norm{\overline{\vek H}_{j}^{(B)}\vek Y - \vek E_{L}^{[(j+1)L]}\vek S_{0}}_{F},\label{eqn.block-GMRES-subprob}
\eea
where $\norm{\cdot}_{F}$ is the Frobenius norm.  
Updates such as $\vek T_{j}^{(G)}$ and $\vek T_{j}^{(F)}$ are often called \emph{corrections} and the
subspaces from which they are drawn are called \emph{correction subspaces}.
There has been a great deal of research
on the convergence properties of block methods such as block GMRES; see, e.g., 
\cite{Gutknecht2007,L.2003,Simoncini.Conv-Block-GMRES}.

In the case $L=1$, block Krylov methods reduce to the well-described single-vector Krylov subspace
methods;
see, e.g., \cite[Section 6.3]{Saad.Iter.Meth.Sparse.2003} and 
\cite{szyld.simoncini.survey.2007}.  In this case, we drop the superscript $(B)$ and write
$\overline{\vek H}_{j} := \overline{\vek H}_{j}^{(B)}$.  The block Arnoldi method simplifies to
a simpler Gramm-Schmidt process in which the block entries $\vek H_{i,j}$ of $\overline{\vek H}_{j}$ reduce
to scalars, now denoted with lower-case $h_{ij}\in\C$.
Then using the scalar version of \eqref{eqn.block-arnoldi-relation}, single-vector FOM can be derived
from \eqref{eqn.block-galerkin-cond} which leads to the formulation
\be\nn
	\vek x_{j}^{(F)} = \vek x_{0} + \vek t_{j}^{(F)}\mwhere \vek t_{j}^{(F)} = \vek V_{j}\vek y_{j}^{(F)}\mand \vek H_{j}\vek y_{j}^{(F)} = \beta\vek e_{1}^{[j]},
\ee
where $\beta = \norm{\vek F_{0}}$ is the $2$-norm of the single-vector residual, and $\vek e_{J}^{[I]}\in\C^{I}$ is
the $J$th Cartesian basis vector in $\C^{I}$.
Similarly for single-vector GMRES, we can use \eqref{eqn.block-arnoldi-relation}, combined with 
\eqref{eqn.block-min-res-cond} to yield the formulation
\be\nn
	\vek x_{j}^{(G)} = \vek x_{0} + \vek t_{j}^{(G)}\mwhere \vek t_{j}^{(G)} = \vek V_{j}\vek y_{j}^{(G)}\mand \vek y_{j}^{(G)}=\argmin{\vek y\in\C^{j}}\norm{\overline{\vek H}_{j}\vek y - \beta\vek e_{1}^{[j+1]}}.
\ee

In the case $L=1$, 
if at some iteration $j$ we have $\CK_{j-1}(\vek A,\vek F_{0}) = \CK_{j}(\vek A,\vek F_{0})$ 
(i.e., $\dim \CK_{j}(\vek A,\vek F_{0}) = j-1 < j$), 
then we have reached an invariant subspace, and both GMRES and FOM will produce
an exact solution at that iteration.  In this case, $j-1$ is called the \emph{grade} of the
pair $\prn{\vek A,\vek F_{0}}$, denoted $\nu(\vek A,\vek F_{0})$.  This notion of grade has been extended to the case 
$L>1$ \cite{GS.2009}; however, the situation is a bit more complicated.  It can
occur that $\dim \K_{j}(\vek A,\vek F_{0}) < jL$ without convergence for all 
right-hand sides (in other
words, without having reached the block grade of $\vek A$ and $\vek F_{0}$,
the iteration at which we reach an invariant subspace).  
It may be
that we have convergence for some or no right-hand sides.  In this case, dependent
block Arnoldi vectors are generated and there must be some procedure in place to 
gracefully handle this situation for reasons of stability.  The dependence of block
Arnoldi vectors and methods for handling this dependence have been discussed extensively
in the literature; see, e.g., 
\cite{B.2000,FM.1997,Gutknecht2007,GS.2009,L.2006,OLeary1980,S.2004,S.2014}, and
general convergence analysis of block methods has been presented in, e.g., 
\cite{L.2003,Simoncini.Conv-Block-GMRES,Simoncini1995}.
In this paper, we consider only the case that dependent basis vectors are replaced with random vectors.

\subsection{Relationships between pairs of projection methods}\label{section.proj-meth-relation}

Pairs of methods such FOM and GMRES which are derived from a Galerkin and minimum
residual projection, respectively, over the same space are closely related.  The analysis
of Brown \cite{B.1991} characterized this relationship in the case of FOM and GMRES when $L=1$.
We state here a theorem encapsulating the results relevant to this work.  First, though,
note that in FOM at iteration $j$, we must solve a linear system involving $\vek H_{j}$.
Thus, if $\vek H_{j}$ is singular, the $j$th FOM iterate does not exist.  We define 
$\widetilde{\vek x}_{j}^{(F)}$ to be the \emph{generalized FOM approximation} through
\be\label{eqn.gen-fom-def}
	\widetilde{\vek x}_{j}^{(F)} = \vek x_{0} + \widetilde{\vek t}_{j}^{(F)}\mwhere \widetilde{\vek t}_{j}^{(F)}=\vek V_{j}\widetilde{\vek y}_{j}^{(F)}\mand \widetilde{\vek y}_{j}^{(F)} = \vek H_{j}^{\dagger}\prn{\beta\vek e_{1}^{[j]}},
\ee
where $\vek H_{j}^{\dagger}$ is the Moore-Penrose pseudoinverse of $\vek H_{j}$.
In the case that $\vek H_{j}$ is nonsingular, we have that $\widetilde{\vek x}_{j}^{(F)} = \vek x_{j}^{(F)}$, but $\widetilde{\vek x}_{j}^{(F)}$ is well-defined in the case that 
$\vek x_{j}^{(F)}$ does not exist.
In this case $\widetilde{\vek y}_{j}^{(F)}$ minimizes 
$\norm{\vek H_{j}\vek y - \beta\vek e_{1}^{[j]}}$
and has minimum norm of all possible minimizers. The following theorem combines two results
proven by Brown in \cite{B.1991}.
\bthm\label{thm.brown-hj-singular}
	The matrix $\vek H_{j}$ is singular (and thus $\vek x_{j}^{(F)}$ does not exist) if and only
	if GMRES stagnates at iteration $j$ with $\vek x_{j}^{(G)} = \vek x_{j-1}^{(G)}$.  
	Furthermore, in the case that $\vek H_{j}$ is singular, we have 
	$\widetilde{\vek x}_{j}^{(F)}=\vek x_{j}^{(G)}$.\footnote{Note that Brown in \cite{B.1991}
	did not use the expression ``generalized FOM approximation''.  
	He calls it the least squares solution and proves it's equivalence to 
	the stagnated $\vek x_{j}^{(G)}$}
\ethm

Thus in the GMRES stagnation case, it is shown that the two methods are ``equivalent'',
if we consider the generalized formulation of FOM.  However, the 
relationship persists in the case that $\vek H_{j}$ is nonsingular as show in, e.g.,
\cite{Saad.Iter.Meth.Sparse.2003}.  In the same text, the following proposition is also shown.
\bprop\label{prop.gmres-single-stag}
	Let $\vek x_{j}^{(G)}$ and $\vek x_{j}^{(F)}$ be the the $j$th GMRES and FOM 
	approximations to the solution of \eqref{eqn.Axb} over the 
	correction subspace $\CK_{j}(\vek A,\vek F_{0})$.  Then we can write
	$\vek x_{j}^{(G)}$ as the following convex combination,
	\be\label{eqn.FOM-GMRES-rel}
		\vek x_{j}^{(G)} = c_{j}^{2}\vek x_{j}^{(F)} + s_{j}^{2}\vek x_{j-1}^{(G)} 
	\ee
	where $s_{j}$ and $c_{j}$ are the $j$th Givens sine and cosine, respectively,
	 obtained from
	annihilating the entry $h_{j+1,j}$ while forming the QR-factorization of
	$\overline{\vek H}_{j}$.
\eprop

One proves this by studying the differences between the QR-factorizations of the 
rectangular $\overline{\vek H}_{j}\in\C^{(j+1)\times j}$ and square 
$\vek H_{j}\in\C^{j\times j}$ generated by the single-vector Arnoldi process.  
The relation \eqref{eqn.FOM-GMRES-rel} reveals information about
GMRES stagnation and its relationship to FOM.  If 
$\vek x_{j}^{(G)} = \vek x_{j-1}^{(G)}$, then we have that $c_{j}=0$
which implies that 	$\vek H_{j}$ is singular and $\vek x_{j}^{(F)}$ does
not exist.  In this case, \eqref{eqn.FOM-GMRES-rel} can be thought of as still valid, 
in the sense that $s_{j}=1$, and \eqref{eqn.FOM-GMRES-rel} reduces to
$\vek x_{j}^{(G)} = \vek x_{j-1}^{(G)}$ if we replace $\vek x_{j}^{(F)}$ with $\widetilde{\vek x}_{j}^{(F)}$.

This characterization of the relationship is not only important for understanding how
these two methods behave at each iteration.  They also reveal that FOM can suffer
from stability issues when GMRES is close to stagnation as the matrix $\vek H_{j}$ is
nearly singular (poorly conditioned) in this case.  Whereas the residual curve of 
GMRES is monotonically nonincreasing, we see spikes in the FOM residual norm corresponding
to periods of near stagnation in the GMRES method.  These so-called ``peaks'' of 
residual norms of FOM and their relation to ``plateaus'' of the residual norms of
GMRES have been previously studied; see, e.g., \cite{CG.1996, C.1994,W.1995,Zw.1994}.
  Of particular 
interest is the observation by Walker \cite{W.1995} 
that the GMRES method can be seen as the 
result of a ``residual smoothing'' of the FOM residual.  Similar observations extend to
other pairings, such as QMR and BiCG. 

\section{Main Results}\label{section.main-results}

When $L>1$, block GMRES and block FOM also fit into the framework of a Galerkin/minimization pairing.  
Thus, it is natural that 
stagnation of block GMRES and behavior of the block FOM algorithm would exhibit the
same relationship, 
using a generalized block FOM iterate defined similar to \eqref{eqn.gen-fom-def}.  
However, this interaction
is more complicated for a block method.  There are interactions between the 
different approximations to individual systems.  
As such, the generalization of stagnation to the block GMRES setting must be done correctly.  We introduce
two definitions.  
\bdefin
At iteration $j$, we call the situation in which $\vek  X_{j}^{(G)} = \vek X_{j-1}^{(G)}$ \emph{total stagnation}.
We call the situation in which some columns
of the block GMRES approximation have stagnated but not all columns \emph{partial stagnation}.
Let $\I$ denote an
	indexing set such that
	$\I\subsetneq\curl{1,2,\ldots, L}$, and let $\overline{\I} = \curl{1,2,\ldots, L}\setminus\I$. 
	For $\vek F\in\C^{J\times L}$, let $\vek F\prn{:,\I}\in\C^{J\times\ab{\I}}$
	have as columns those from $\vek F$ corresponding to indices in $\I$.
	Then partial stagnation refers to the situation in which we have 
	\be\label{eqn.partial-stag}
		\vek X_{j}^{(G)}\prn{:,\I} = \vek X_{j-1}^{(G)}\prn{:,\I}\mbut\vek X_{j}^{(G)}\prn{:,i} \neq\vek X_{j-1}^{(G)}\prn{:,i}\mforeach i\in\overline{\I}.
	\ee
	\edefin
	Total stagnation is analogous to stagnation of GMRES in the single-vector case, as 
characterized in \cite{B.1991}, but partial stagnation has no single-vector analog.
	Both total and partial stagnation can occur for multiple reasons.  Total block GMRES stagnation can occur when block 
	GMRES has converged, i.e., $\vek X_{j}^{(G)} = \vek X$,  implying (if $j$ is the first iteration for which this occurs) 
	from \cite[Theorem 9]{GS.2009}, that we have that $j=\nu(\vek A, \vek F_{0})$ and 
	$\dim \K_{j+k}(\vek A,\vek F_{0})=\dim \K_{j}(\vek A,\vek F_{0})$ for all $k>0$.   This case is trivial
	and will not be considered.
	If there is no breakdown of the block Arnoldi process (the rank of the block residual is $L$), then an occurrence 
	of total stagnation is the
	block analog of single-vector GMRES stagnation.  We prove in this case that Theorem \ref{thm.brown-hj-singular}
	has a block analog; c.f., Corollary \ref{cor.Hsing-GMRES-total-stag} and Corollary \ref{cor.bl-GMRES-stag-FOM-stag}.
	
	Partial stagnation has no direct analog to the single-vector case. Partial stagnation can occur when for 
	column $i$, the system is exactly solved with \linebreak$\vek X_{j}^{(G)}(:,i) = \vek X(:,i)$.  This implies that 
	$\vek F_{0}(:,i) - \vek A\vek W_{j}\vek Y_{j}^{(G)}(:,i)=0$, which implies that 
	\be\nn
		\dim \prn{\CR(\vek F_{0}) \cap \vek A\K_{j}(\vek A, \vek F_{0})}=1
	\ee	
	(see, e.g., \cite{RS.2006}) and that a dependent
	Arnoldi vector has been produced.  In this case, one can treat this with one of the referenced 
	strategies; see, e.g., \cite{OLeary1980,B.2000,ABFH.2000,BF.2013,FM.1997,SO.2010}.
	
	This is a specific instance of block Arnoldi process breakdown. At iteration $j$, the process breaks down when
	the matrix $\bbmat \vek B & \vek A\vek B & \cdots & \vek A^{j-1}\vek B \ebmat$ is rank deficient which is equivalent
	to saying $\dim\prn{\CR(\vek X)\cap \K_{j}(\vek A,\vek F_{0})} = \dim\calN(\vek R_{j}) > 0$. 
	In this case, $\K_{j}\prn{\vek A,\vek F_{0}}$ contains a linear combination of the columns of $\vek X$ 
	\cite{Nikishin1995,RS.2006}.  
	It has also been observed \cite{RS.2006} that a dependent Arnoldi vector can be generated without 
	the convergence of any of the columns.  
	
	In the case that there has been no breakdown of
	the block Arnoldi process we show that partial stagnation is actually a special case of a more general situation in which
	a part of the Krylov subspace does not contribute to the GMRES minimization process and the dimension of this
	subspace corresponds to the dimension of the null space of the rank-deficient FOM matrix $\vek H^{(B)}_{j}$;
	c.f., Theorem \ref{thm.H-sing} and Theorem \ref{thm.r-dim-update-bl-GMRES-FOM} below.
	
		We derive a relationship for block GMRES and block FOM which is a 
	generalization of \eqref{eqn.FOM-GMRES-rel} and is valid even in the case that
	$\vek H_{j}^{(B)}$ is singular.  Thus, as in \eqref{eqn.gen-fom-def}, we generalize
	the definition of the block FOM approximation to be compatible with a singular 
	$\vek H_{j}^{(B)}$, i.e., 
	\bea
		\widetilde{\vek X}_{j}^{(F)} = \vek X_{0} + \widetilde{\vek T}_{j}^{(F)}& \mwhere & \widetilde{\vek T}_{j}^{(F)}=\vek W_{j}\widetilde{\vek Y}_{j}^{(F)}\nn\\ &\mand & \widetilde{\vek Y}_{j}^{(F)} = \prn{\vek H_{j}^{(B)}}^{\dagger}\prn{\vek E_{L}^{[jL]}\vek S_{0}}\label{eqn.block-FOM-generalized}
	\eea
	where $\prn{\vek H_{j}^{(B)}}^{\dagger}$ is the Moore-Penrose pseudoinverse 
	of $\vek H_{j}^{(B)}$.
In the case that $\vek H_{j}^{(B)}$ is nonsingular, we have that $\widetilde{\vek X}_{j}^{(F)} = \vek X_{j}^{(F)}$, but $\widetilde{\vek X}_{j}^{(F)}$ is well-defined in the case that 
$\vek X_{j}^{(F)}$ does not exist.
In this case $\widetilde{\vek Y}_{j}^{(F)}$ minimizes 
$\norm{\vek H_{j}^{(B)}\vek Y - \vek E_{L}^{[jL]}\vek S_{0}}_{F}$
and has minimum norm of all possible minimizers. 
	As in 
	\eqref{eqn.gen-fom-def}, this definition reduces to the standard formulation of 
	the FOM approximation in the case that $\vek H_{j}^{(B)}$ is nonsingular. 
	In the single-vector case,
	to prove \cite[Lemma 6.1]{Saad.Iter.Meth.Sparse.2003}, expressions are derived for the
	inverses of upper-triangular matrices.  We need to obtain similar identities here.
	However, we want our derivation to be compatible with the case that $\vek H_{j}^{(B)}$
	is singular.

To characterize both types of stagnation requires us to follow the work in \cite{B.1991},
generalizing to the block Krylov subspace case.  We also need to generalize 
\eqref{eqn.FOM-GMRES-rel} to the block GMRES/FOM setting. This is quite useful 
in extending the work in \cite{B.1991} and also of general interest.  

\subsection{GMRES and FOM from a particular perspective}\label{section.other-perspective}
We discuss briefly the known results for the relationship of single-vector GMRES and (generalized) FOM. This discussion closely 
relates to the discussion and results on ascent directions in, e.g., \cite{B.1991}.
It has been shown that at the $j$th iteration the approximations $\vek x_{j}^{(G)}$ and $\widetilde{\vek x}_{j}^{(F)}$ can both
be related to the $(j-1)st$, with
\bea
	\vek x_{j}^{(G)} = \vek x_{j-1}^{(G)} + \vek s_{j}^{(G)}& \mand &\widetilde{\vek x}_{j}^{(F)} = \vek x_{j-1}^{(G)} + \widetilde{\vek s}_{j}^{(F)}\label{eqn.GMRES-FOM-iter-update}\\\mwhere \vek s_{j}^{(G)}=\vek V_{j}\vek y^{(G)}_{\vek s_{j}}\in\CK_{j}(\vek A, \vek F_{0})&\mand &\widetilde{\vek s}_{j}^{(F)}=\vek V_{j}\widetilde{\vek y}_{\vek s_{j}}^{(F)}\in\CK_{j}(\vek A, \vek F_{0})\nn
\eea
where $\widetilde{\vek y}_{\vek s_{j}}^{(F)}$ and $\vek y^{(G)}_{\vek s_{j}}$ are representations of the 
generalized FOM and GMRES \emph{progressive} corrections from $\CK_{j}(\vek A, \vek F_{0})$.
The next proposition follows directly.
\bprop\label{prop.GF-resid-aug}
	The GMRES and generalized FOM updates $\vek y^{(G)}_{\vek s_{j}}$ and $\widetilde{\vek y}_{\vek s_{j}}^{(F)}$ respectively satisfy the minimizations
	\bea
		\vek y^{(G)}_{\vek s_{j}} &=& \argmin{\vek y\in\Cn}\norm{\bbmat \beta\vek e_{1}^{[j]} - \overline{\vek H}_{j-1}\vek y_{j-1}^{(G)} \\ 0 \ebmat - \overline{\vek H}_{j}\vek y} \mand \label{eqn.GMRES-iter-argmin}\\  \widetilde{\vek y}_{\vek s_{j}}^{(F)} &=& \argmin{\vek y\in\Cn}\norm{\beta\vek e_{1}^{[j]} - \overline{\vek H}_{j-1}\vek y_{j-1}^{(G)} - \vek H_{j}\vek y}.
		\label{eqn.FOM-iter-argmin}
	\eea
\eprop

\bproof
	To prove \eqref{eqn.GMRES-iter-argmin}, one simply inserts the expression 
	for $\vek x_{j}^{(G)}$ from \eqref{eqn.GMRES-FOM-iter-update} into the residual and applies the GMRES Petrov-Galerkin condition 
	\eqref{eqn.block-min-res-cond}.  To prove \eqref{eqn.FOM-iter-argmin}, one begins similarly, by substituting the expression 
	for $\vek x_{j}^{(F)}$ from \eqref{eqn.GMRES-FOM-iter-update} into the residual and applying the FOM Galerkin condition 
	\eqref{eqn.block-galerkin-cond}.  In this case, if $\vek H_{j}$ is nonsingular, this is equivalent to solving the linear system
	\be\label{eqn.FOM-aug-sys}
		\vek H_{j}\widetilde{\vek y}_{\vek s_{j}}^{(F)} = \beta\vek e_{1}^{[j]} - \overline{\vek H}_{j-1}\vek y_{j}^{(G)}.
	\ee
	In the case that $\vek H_{j}$ is singular (the $j$th FOM approximation does not exist), we set 
	\be\label{eqn.FOM-aug-LS}
		\widetilde{\vek y}_{\vek s_{j}}^{(F)} = \vek H_{j}^{\dagger}\prn{\beta\vek e_{1}^{[j]} - \overline{\vek H}_{j-1}\vek y_{j}^{(G)}}.
	\ee
	In either case, we have that $\widetilde{\vek y}_{\vek s_{j}}^{(F)}$ is the minimizer of \eqref{eqn.FOM-iter-argmin}, yielding
	the result.
\eproof
The result on FOM is \cite[Theorem 3.3]{B.1991} but stated differently.  
This
formulation allows us to discuss the GMRES and FOM 
at iteration $j$ using the $(j-1)$st GMRES minimization.  We see that
the GMRES method least-squares problem simply grows by one dimension when we go from iteration $j-1$ to $j$.
However, at iteration $j$, imposing the FOM Galerkin condition 
\eqref{eqn.block-galerkin-cond} is equivalent to an augmentation of the $(j-1)$st GMRES least squares matrix. 
 This augmented matrix is square.  If it is nonsingular, then the $j$th FOM approximation exists and
 we solve the augmented system \eqref{eqn.FOM-aug-sys}.  If the augmented matrix is singular, then the generalized FOM approximation
 is computed by solving the least squares problem \eqref{eqn.FOM-aug-LS}.  In the case of single-vector GMRES and FOM, this
is not necessary to characterize their relationship. However, in the case of block GMRES and block FOM, we can better discuss
a generalization to the more complicated block Krylov subspace situation.  

\subsection{The QR-Factorization of the block upper Hessenberg matrices}\label{section.bl-hess-QR}
We begin by describing the structure of the 
QR-factorizations of the square and rectangular block Hessenberg matrices.
\blem
	Let $\overline{\vek R}_{j}\in\C^{j+1\times j}$ and 
	$\widehat{\vek R}_{j}\in\C^{j\times j}$
	be the R-factors of the respective QR-factorizations of $\overline{\vek H}_{j}^{(B)}$ 
	and $\vek H_{j}^{(B)}$, and let $\vek R_{j}$ be the $j\times j$ non-zero block of
	$\overline{\vek R}_{j}$.  
	Then $\overline{\vek R}_{j}$ and $\widehat{\vek R}_{j}$ both have as their
	upper left $j-1\times j-1$ block the R-factor of
	the QR-factorization of $\overline{\vek H}_{j-1}^{(B)}$, i.e., $\vek R_{j-1}$.  
	Furthermore, the structures of 
	$\overline{\vek R}_{j}$ and $\widehat{\vek R}_{j}$, respectively, are,
	\be\label{eqn.R-structures}
		\overline{\vek R}_{j}=\bbmat \vek R_{j-1}  & \vek Z_{j}\\ 
														                   & \vek N_{j} \\
														                   & 
										  \ebmat\mand
						\widehat{\vek R}_{j}= \bbmat\vek R_{j-1} & \vek Z_{j}\\ 
														        & \widehat{\vek N}_{j}
										   \ebmat,
	\ee
	where $\vek Z_{j}\in\C^{(j-1)L\times L}$ and
	$\vek N_{j},\widehat{\vek N}_{j}\in\C^{L\times L}$ are upper triangular.
\elem
\bproof
	Let  
	$\vek Q_{i}^{(j+1)}  \in\C^{j+1\times j+1}$ be orthonormal transformation which annihilates all subdiagonal
	entries in columns $i-1+1$ to $i$ of $\overline{\vek H}_{j}^{(B)}$ and effects no other rows
	so that we can write 
	\small
	\be\nn
	\vek Q_{j-1}^{(j+1)}\cdots\vek Q_{1}^{(j+1)}\overline{\vek H}_{j}^{(B)} = \bbmat\vek R_{j-1} & \vek Z_{j}\\     & \widehat{\vek H}_{j,j}\\ & \vek H_{j+1,j}\ebmat \mand\vek Q_{j-1}^{(j)}\cdots\vek Q_{1}^{(j)}\vek H_{j}^{(B)} = \bbmat\vek R_{j-1} & \vek Z_{j} \\  & \widehat{\vek H}_{j,j}\ebmat.
	\ee
	\normalsize
	Let $\widehat{\vek Q}_{j}^{(j)}\in\C^{j\times j}$ be the orthogonal transformation which annihilate the lower subdiagonal entries the
	block $\widehat{\vek H}_{j,j}$ in $\vek Q_{j-1}^{(j)}\cdots\vek Q_{1}^{(j)}\vek H_{j}^{(B)}$ and effects no other rows.
	Then we have 
	\be\label{eqn.Rj-orthog-trans}
		\overline{\vek R}_{j} = \vek Q_{j}^{(j+1)}\bbmat\vek R_{j-1} & \vek Z_{j}\\     & \widehat{\vek H}_{j,j}\\ & \vek H_{j+1,j}\ebmat\mand \widehat{\vek R}_{j} = \widehat{\vek Q}_{j}^{(j)}\bbmat\vek R_{j-1} & \vek Z_{j} \\  & \widehat{\vek H}_{j,j}\ebmat,
	\ee 
	and the Lemma is proven.
\eproof
Thus, the two core problems which must be solved at every iteration of block GMRES
and block FOM can be written
\be\label{eqn.Y-gmres-solve}
\bbmat \vek R_{j-1}  & \vek Z_{j}\\ 
								  & \vek N_{j}\ebmat\vek Y_{j}^{(G)} = (\vek Q_{j}^{(j+1)}\cdots\vek Q_{1}^{(j+1)}\vek E_{j+1}^{[j+1]}\vek S_{0})_{1:j} 
\ee
\be\label{eqn.Y-fom-solve}
\mand\bbmat\vek R_{j-1} & \vek Z_{j}\\ 
														        & \widehat{\vek N}_{j}
										   \ebmat\vek Y_{j}^{(F)} = \widehat{\vek Q}_{j}^{(j)}\vek Q_{j-1}^{(j)}\cdots\vek Q_{1}^{(j)}\vek E_{j}^{[j]}\vek S_{0}.
\ee
It is also straightforward to show that the block right-hand sides of these core problems
are related. If
	\be\nn
		\vek G_{j}^{(G)} =  (\vek Q_{j}^{(j+1)}\cdots\vek Q_{1}^{(j+1)}\vek E_{L}^{[j+1]}\vek S_{0})_{1:j} \mand \vek G_{j}^{(F)} = \widehat{\vek Q}_{j}^{(j)}\vek Q_{j-1}^{(j)}\cdots\vek Q_{1}^{(j)}\vek E_{L}^{[j+1]}\vek S_{0},
	\ee
	then $\vek G_{j}^{(G)}$ and $\vek G_{j}^{(F)}$ are equal for the first $j-1$
	rows, with
	\be\label{eqn.G-expressions}
		\vek G_{j}^{(G)} = \bbmat\vek G_{j-1}^{(G)}\\ \vek C_{j}\ebmat\mand\vek G_{j}^{(F)}=\bbmat\vek G_{j-1}^{(G)}\\ \widehat{\vek C}_{j}\ebmat
	\ee
	where we have that
	\be\label{eqn.GjG-structure}
		\vek G_{j}^{(G)} = \prn{\vek Q_{j}^{(j+1)}\bbmat\vek G_{j-1}^{(G)}\\ \widetilde{\vek C}_{j}\\ \vek 0\ebmat}_{1:j} = \prn{\bbmat\vek G_{j-1}^{(G)}\\ {\vek C}_{j}\\ \boldsymbol\ast\ebmat}_{1:j}.
	\ee
This is a consequence of the structure of the orthogonal transformations used to define these vectors.
	It is important to pause here for a moment to discuss the $L\times L$ matrices
	$\vek C_{j}$, $\widehat{\vek C}_{j}$, and $\widetilde{\vek C}_{j}$ and characterize
	if and when they are full rank.  At times for convenience, we refer to these matrices as the ``$\vek C$-matrices''.
	\blem\label{lemma.Ctilde}
		We have that $\rank \widetilde{\vek  C}_{j} = \rank \vek F_{j-1}^{(G)}$;
		and, in particular,  if \linebreak$\dim \K_{j-1}(\vek A,\vek F_{0})=(j-1)L$,
		we have that, $\widetilde{\vek  C}_{j}$ is nonsingular.  
	\elem
	\bproof
		Let $\vek Y_{j-1}^{(G)}$ be the solution to the block GMRES least squares 
		subproblem \eqref{eqn.block-GMRES-subprob} but for iteration $j-1$.  Let
		\be\nn
			\vek F_{j-1}^{(G)} = \vek B - \vek A\vek X_{j-1}^{(G)} = \vek W_{j}\prn{\overline{\vek H}_{j-1}\vek Y_{j-1}^{(G)} - \vek E_{L}^{[j]}\vek S_{0}}.
		\ee
		By assumption \eqref{eqn.Y-gmres-solve} has a solution at iteration $j-1$,
		and thus 
		\be\nn
			\overline{\vek H}_{j-1}\vek Y_{j-1}^{(G)} - \vek E_{L}^{[j]}\vek S_{0} = \vek Q_{j-1}^{\ast}\bbmat\vek G_{j-1}^{(G)}\\ \vek C_{j}\\ \vek 0\ebmat - \vek E_{L}^{[j]}\vek S_{0}.
		\ee
		where $\overline{\vek Q}_{j-1} = \vek Q_{j-1}^{(j)}\cdots\vek Q_{1}^{(j)}$. 
		Since $\vek W_{j-1}$ and $\overline{\vek Q}_{j-1}$ are both full rank, we have
		\bea
			\rank\vek F_{j-1}^{(G)} = \rank \vek Q_{j-1}\vek W_{j}^{\ast}\vek F_{j}^{(G)} = \rank\prn{\bbmat\vek G_{j-1}^{(G)}\\ \vek C_{j-1}\\ \vek 0\ebmat - \bbmat\vek G_{j-1}^{(G)}\\ \vek C_{j-1}\\ \widetilde{\vek C}_{j}\ebmat} &=& \rank\bbmat\vek 0\\\vek 0 \\ -\widetilde{\vek C}_{j}\ebmat\nn \\ &=& \rank\widetilde{\vek C}_{j}.\nn
		\eea
		If we assume that the block Arnoldi method has not produced any dependent
		basis vectors, then we know from \cite[Section 2, Corollary 1]{RS.2006} that $\vek F_{j-1}^{(G)}$ is full-rank
		 meaning $\widetilde{\vek C}_{j}$ is  nonsingular.
	\eproof
	
	From this, we can similarly characterize the ranks of
	${\vek C}_{j}$ and $\widehat{\vek C}_{j}$ which are closely related to 
	$\widetilde{\vek C}_{j}$.
	\blem\label{lemma.Chat}
		We have that $\rank\widehat{\vek C}_{j}=\rank\widetilde{\vek C}_{j}$.  In 
		particular, if $\dim\K_{j-1}(\vek A,\vek F_{0}) = (j-1)L$, we have that $\widehat{\vek C}_{j}$ is square
		and nonsingular.
	\elem
	\bproof
		Let $\widehat{\vek Q}_{j}^{(j)} = \bbmat\vek I_{j-1} & \\ & \widehat{\vek Q}_{j}^{(b)}\ebmat$ where $\widehat{\vek Q}_{j}^{(b)}\in\C^{L\times L}$ is the orthogonal 
		transformation  such that the second equation of \eqref{eqn.Rj-orthog-trans} holds.  
		Then from \eqref{eqn.GjG-structure}, we have 
		$\widehat{\vek C}_{j} = \widehat{\vek Q}_{j}^{(b)}\widetilde{\vek C}_{j}$. If $\widetilde{\vek C}_{j}$ has 
		full rank, the second statement follows.
	\eproof
	
	We can prove a similar result for $\vek C_{j}$, which will be used later to verify the nonsingularity of
	$\vek C_{j}$ under certain conditions.
	\blem\label{lemma.C}
		Let 
		\be\nn
			\vek Q_{j}^{(j+1)} = \bbmat\vek I_{\CJ(j-2)} &  &\\ & \vek Q_{j}^{(11)} & \vek Q_{j}^{(11)} \\ & \vek Q_{j}^{(21)} & \vek Q_{j}^{(22)} \ebmat
		\ee		
		with $\vek Q_{j}^{(11)}\in\C^{L\times L}$, $\vek Q_{j}^{(12)}\in\C^{L\times L}$, 
		$\vek Q_{j}^{(21)}\in\C^{L\times L}$, and $\vek Q_{j}^{(22)}\in\C^{L\times L}$.  
		In general, we have $rank\,\vek C_{j}\leq \min \curl{\rank\vek Q_{j}^{(11)},\rank\widetilde{\vek C}_{j}}$.
		If  $\dim\K_{j-1}(\vek A,\vek F_{0})~=~(j~-~1)L$, we have 
		$\vek C_{j}$ is singular if and only if $\vek Q_{j}^{(11)}$ is singular.
	\elem
	\bproof
		From \eqref{eqn.GjG-structure} we have that $\vek C_{j} = \vek Q_{j}^{(11)}\widetilde{\vek C}_{j}$.  
		The general result comes from basic inequality results for ranks of products of matrices; see, e.g., \cite[Chapter 0]{Horn1985}.
		If we assume  $\dim\K_{j-1}(\vek A,\vek F_{0})~=~(j~-~1)L$, then we know that $\widetilde{\vek C}_{j}$ has full rank,
		and 
		the second result (in both directions) follows.
	\eproof
	
	We see that the ranks of $\widetilde{\vek C}_{j}$ and $\widehat{\vek C}_{j}$ are directly connected to block Arnoldi
	breakdown at iteration $j-1$.  Later in Section \ref{section.no-breakdown}, we assume no breakdown, thus both $\widetilde{\vek C}_{j}$ and $\widehat{\vek C}_{j}$ are nonsingular.  In Section \ref{section.breakdown-case},
	we assume that the block Arnoldi process produces dependent vectors at iteration $j$ which are replaced with random vectors.  Thus, at iteration $j$, 
	both $\widetilde{\vek C}_{j}$ and $\widehat{\vek C}_{j}$ are still nonsingular, and their dimensions do not change at subsequent iterations.
		
	We now turn to solving \eqref{eqn.Y-fom-solve} and either
	solving \eqref{eqn.Y-gmres-solve} or obtaining the generalized least squares solution
	if $\widehat{\vek R}_{j}$ is singular.
	 Since $\vek R_{j}$ is
	nonsingular, we simply compute the actual inverse while  for $\widehat{\vek R}_{j}$,
	we compute the pseudo-inverse.
	These are both straightforward
	generalizations of the identities used in the proof of 
	\cite[Lemma 6.1]{Saad.Iter.Meth.Sparse.2003}, though verifying the structure
	of the Moore-Penrose
	pseudo-inverse identity requires a bit of thought.  Let us recall briefly the following
	definition which can be found in, e.g., \cite[Section 2.2]{EHN.1996-book},
	\bdefin\label{def.Moore-Penrose}
		Let $T:\CX\rightarrow\CY$ be a bounded linear operator between Hilbert spaces.
		Let $\calN(T)$ denote the null space and $\CR(T)$ denote the range of $T$ and
		define $\widetilde{T}:\calN(T)^{\perp}\rightarrow\CR(T)$ to be the invertible 
		operator such that $\widetilde{T}x = Tx$ for all $x~\in~\calN(T)^{\perp}$.
		Then we call the operator $T^{\dagger}$ the Moore-Penrose pseudo-inverse if
		it is the unique operator satisfying
		\begin{enumerate}
			\item\label{item.MP-cond1} $\restr{T^{\dagger}}{\CR(T)} = \widetilde{T}^{-1}$
			\item\label{item.MP-cond2} $\restr{T^{\dagger}}{\CR(T)^{\perp}} = 0_{op}$
		\end{enumerate}
		where $0_{op}$ is the zero operator.
	\edefin
	
	This definition is more general than the matrix-specific definition given in, e.g., \cite[Section 5.5.2]{GV.2013}.
	We choose to follow Definition \ref{def.Moore-Penrose} as it renders the proof of the following lemma 
	less dependent on many lines of block matrix calculations, but of course the theoretical results are the same.
	\blem\label{lem.R-inverse}
		The inverse and pseudo-inverse, respectively, of $\vek R_{j}$ and $\widehat{\vek R}_{j}$ 
		can be directly constructed from the identities \eqref{eqn.R-structures}, i.e.,
		\be\label{eqn.R-inverse}
			\vek R_{j}^{-1} = \bbmat\vek R_{j-1}^{-1}& -\vek R_{j-1}^{-1}\vek Z_{j}\vek N_{j}^{-1}\\ & \vek N_{j}^{-1}\ebmat\mand\widehat{\vek R}_{j}^{\dagger}=\bbmat\vek R_{j-1}^{-1}& -\vek R_{j-1}^{-1}\vek Z_{j}\widehat{\vek N}_{j}^{\dagger}\\ & \widehat{\vek N}_{j}^{\dagger}\ebmat,
		\ee
		where $\widehat{\vek N}_{j}^{\dagger}$ is the Moore-Penrose pseudo-inverse of
		$\widehat{\vek N}_{j}$.
	\elem
	\bproof
		The expression for $\vek R_{j}^{-1}$ can be directly verified by left and right multiplication.  To verify
		the expression for $\widehat{\vek R}_{j}^{\dagger}$, we must verify the two conditions listed in 
		Definition \ref{def.Moore-Penrose}.  

		To verify Condition \ref{item.MP-cond1}, we first construct a basis for $\calN(\widehat{\vek R}_{j})$.
		Observe that under our assumption that $\vek R_{j-1}$ is nonsingular, we have that 
		\be\nn
			\dim\calN(\widehat{\vek R}_{j}) = \dim\calN(\widehat{\vek N}_{j}) = L-r
		\ee		 
		 where $r = \rank(\widehat{\vek N}_{j})$. 
		Let $\curl{\vek y_{i}}_{i=1}^{r}$ be a basis for $\calN(\widehat{\vek N}_{j})^{\perp}$.  Furthermore,
		let $\curl{\vek m_{i}}_{i=1}^{j-1}$ be a basis for $\R^{j-1}$.  Then it follows that 
		\be\nn
			\curl{\bbmat \vek R_{j-1}^{-1}\vek m_{1}\\ \vek 0 \ebmat, \ldots, \bbmat \vek R_{j-1}^{-1}\vek m_{j-1}\\ \vek 0 \ebmat, \bbmat -\vek R_{j-1}^{-1}\vek Z_{j}\vek y_{1}\\ \vek y_{1}\ebmat, \ldots, \bbmat -\vek R_{j-1}^{-1}\vek Z_{j}\vek y_{r}\\ \vek y_{r} \ebmat}
		\ee
		is a basis for $\calN(\widehat{\vek R}_{j})^{\perp}$.  For any 
		$\widehat{\vek x}\in\calN(\widehat{\vek R}_{j})^{\perp}$, we can write 
		\be\nn
		\widehat{\vek x} = \sum_{i=1}^{j-1}\alpha_{i}\bbmat \vek R_{j-1}^{-1}\vek m_{i}\\ \vek 0 \ebmat + \sum_{i=1}^{r}\beta_{i}\bbmat -\vek R_{j-1}^{-1}\vek Z_{j}\vek y_{i}\\ \vek y_{i} \ebmat.
		\ee
		By direct calculation, we see that 
		\be\nn
			\widehat{\vek R}_{j}\widehat{\vek x} = \sum_{i=1}^{j-1}\alpha_{i}\bbmat \vek m_{i}\\\vek 0 \ebmat + \sum_{i=1}^{r}\beta_{i}\bbmat \vek 0\\ \widehat{\vek N}_{j}\vek y_{j} \ebmat,
		\ee
		and applying our prospective pseudo-inverse yields
		\be\nn
			\widehat{\vek R}_{j}^{\dagger}\widehat{\vek R}_{j}\widehat{\vek x} = \sum_{i=1}^{j-1}\alpha_{i}\bbmat \vek R_{j-1}^{-1}\vek m_{i}\\\vek 0 \ebmat + \sum_{i=1}^{r}\beta_{i}\bbmat -\vek R_{j-1}^{-1}\vek Z_{j}\widehat{\vek N}_{j}^{\dagger}\widehat{\vek N}_{j}\vek y_{i}\\ \widehat{\vek N}_{j}^{\dagger}\widehat{\vek N}_{j}\vek y_{i} \ebmat.
		\ee
		Finally, we observe that since $\curl{\vek y_{i}}_{i=1}^{r}$ is a basis for 
		$\calN(\widehat{\vek N}_{j})^{\perp}$,
		we have from Definition~\ref{def.Moore-Penrose} 
		that $\widehat{\vek N}_{j}^{\dagger}\widehat{\vek N}_{j}\vek y_{i}~=~\vek y_{i}$ 
		for all $i$, and 
		thus 
		$\widehat{\vek R}_{j}^{\dagger}\widehat{\vek R}_{j}\widehat{\vek x}  = \widehat{\vek x} $, verifying 
		Condition \ref{item.MP-cond1}.
		
		To verify Condition \ref{item.MP-cond2}, we first observe that 
		\be\nn
			\curl{\bbmat \vek m_{1}\\ \vek 0 \ebmat, \ldots, \bbmat \vek m_{j-1}\\ \vek 0 \ebmat, \bbmat \vek 0\\ \widehat{\vek N}_{j}\vek y_{1}\ebmat, \ldots, \bbmat \vek 0\\ \widehat{\vek N}_{j}\vek y_{r} \ebmat}
		\ee
		is a basis for $\CR(\widehat{\vek R}_{j})$.  Let $\curl{\vek c_{i}}_{i=1}^{L-r}$ 
		be a basis for $\CR(\widehat{\vek N}_{j})^{\perp}$.
		Then it follows that $\curl{\bbmat \vek 0\\ \vek c_{i}\ebmat}_{i=1}^{L-r}$ is a 
		basis for $\CR(\widehat{\vek R}_{j})^{\perp}$.
		Let $\tilde{\vek y} = \sum_{i=1}^{L-r}\gamma_{i}\bbmat \vek 0\\ \vek c_{i}\ebmat$ be an element of 
		$\CR(\widehat{\vek R}_{j})^{\perp}$.  Then we have 
		\be\nn
			\widehat{\vek R}_{j}^{\dagger}\tilde{\vek y} = \sum_{i=1}^{L-r}\gamma_{i}\bbmat -\vek R_{j-1}^{-1}\vek Z_{j}\widehat{\vek N}_{j}^{\dagger}\vek c_{i}\\ \widehat{\vek N}_{j}^{\dagger}\vek c_{i}\ebmat.
		\ee
		It follows directly from Definition \eqref{def.Moore-Penrose} that 
		$\widehat{\vek N}_{j}^{\dagger}\vek c_{i}=0$ for all $i$, and this proves Condition \ref{item.MP-cond2},
		thus proving the the lemma.
	\eproof
	The following corollary technically follows from Lemma \ref{lem.R-inverse}, though it can easily be 
	proven directly.
	\bcor
		If $\vek H_{j}^{(B)}$ is nonsingular, then it follows that $\widehat{\vek R}_{j}^{-1}$ and can be
		written
		\be\nn
			\widehat{\vek R}_{j}^{-1}=\bbmat\vek R_{j-1}^{-1}& -\vek R_{j-1}^{-1}\vek Z_{j}\widehat{\vek N}_{j}^{-1}\\ & \widehat{\vek N}_{j}^{-1}\ebmat.
		\ee
	\ecor
	
	Now we have all the pieces we need to analyze the relationship
	between the block GMRES and block FOM approximations, and we can then
	discuss the implications with respect to stagnation.  
	\subsection{The case of a breakdown-free block Arnoldi process}\label{section.no-breakdown}
	
	We begin this section by discussing block GMRES and block FOM from the same perspective as advocated in
	Section \ref{section.other-perspective}.  
	We have the block analog of Proposition \ref{prop.GF-resid-aug}, and in this case
	we explicitly construct the block analogs of $\vek s_{j}^{(G)}$ and $\vek s_{j}^{(F)}$.
	\blem\label{lemma.block-GF-resid-aug}
		Let $\vek S_{j}^{(G)} = \vek W_{j}\vek Y_{\vek S_{j}}^{(G)}$ and 
		$\widetilde{\vek S}_{j}^{(F)} = \vek W_{j}\widetilde{\vek Y}_{\vek S_{j}}^{(F)}$ both be in $\C^{n\times L}$ such that
		they satisfy the block GMRES and FOM progressive update formulas
		\bea
			\vek X_{j}^{(G)} = \vek X_{j-1}^{(G)} + \vek S_{j}^{(G)}& \mand &\widetilde{\vek X}_{j}^{(F)} = \vek X_{j-1}^{(G)} + \widetilde{\vek S}_{j}^{(F)}\label{eqn.bl-GMRES-FOM-iter-update}.
		\eea
		Then we can write 
		\be\label{eqn.block-Y-update}
			\vek Y_{\vek S_{j}}^{(G)} = \bbmat  -\vek R_{j-1}^{-1}\vek Z_{j}\vek N_{j}^{-1}\vek C_{j} \\ \vek N_{j}^{-1}\vek C_{j} \ebmat\mand \widetilde{\vek Y}_{\vek S_{j}}^{(F)} = \bbmat  -\vek R_{j-1}^{-1}\vek Z_{j}\widehat{\vek N}_{j}^{\dagger}\widehat{\vek C}_{j} \\ \widehat{\vek N}_{j}^{\dagger}\widehat{\vek C}_{j} \ebmat,
		\ee
		and these vectors minimize the two residual
		update equations
		\bea
			\vek Y^{(G)}_{\vek S_{j}} &=& \argmin{\vek Y\in\Cn}\norm{\bbmat \vek E_{1}^{[jL]}\vek S_{0} - \overline{\vek H}^{(B)}_{j-1}\vek Y_{j-1}^{(G)} \\ \vek 0 \ebmat - \overline{\vek H}^{(B)}_{j}\vek Y} \mand \label{eqn.bl-GMRES-iter-argmin}\\  \widetilde{\vek Y}_{\vek S_{j}}^{(F)} &=& \argmin{\vek Y\in\Cn}\norm{\vek E_{1}^{[jL]}\vek S_{0} - \overline{\vek H}^{(B)}_{j-1}\vek Y_{j-1}^{(G)} - \vek H^{(B)}_{j}\vek Y}.
			\label{eqn.bl-FOM-iter-argmin}
		\eea
	\elem
	\bproof
		Combining \eqref{eqn.G-expressions} and \eqref{eqn.R-inverse} to solve 
		\eqref{eqn.Y-gmres-solve} and \eqref{eqn.block-FOM-generalized}
		we have the following expressions for $\vek Y_{j}^{(G)}$ and 
		$\widetilde{\vek Y}_{j}^{(F)}$,
		\small
		\be\nn
			\vek Y_{j}^{(G)} = \bbmat \vek R_{j-1}^{-1}\vek G_{j-1}^{(G)} -\vek R_{j-1}^{-1}\vek Z_{j}\vek N_{j}^{-1}\vek C_{j} \\ \vek N_{j}^{-1}\vek C_{j} \ebmat\mand\widetilde{\vek Y}_{j}^{(F)} = \bbmat \vek R_{j-1}^{-1}\vek G_{j-1}^{(G)} -\vek R_{j-1}^{-1}\vek Z_{j}\widehat{\vek N}_{j}^{\dagger}\widehat{\vek C}_{j} \\ \widehat{\vek N}_{j}^{\dagger}\widehat{\vek C}_{j} \ebmat
		\ee
		\normalsize
		As it can be appreciated, 
		$\vek R_{j-1}^{-1}\vek G_{j-1}^{(G)} = \vek Y_{j-1}^{(G)}$, and it follows that
		\be\nn
			\vek Y_{j}^{(G)} = \bbmat\vek Y_{j-1}^{(G)}\\ \vek 0\ebmat + \bbmat  -\vek R_{j-1}^{-1}\vek Z_{j}\vek N_{j}^{-1}\vek C_{j} \\ \vek N_{j}^{-1}\vek C_{j} \ebmat\mand\widetilde{\vek Y}_{j}^{(F)} = \bbmat\vek Y_{j-1}^{(G)}\\ \vek 0\ebmat +\bbmat  -\vek R_{j-1}^{-1}\vek Z_{j}\widehat{\vek N}_{j}^{\dagger}\widehat{\vek C}_{j} \\ \widehat{\vek N}_{j}^{\dagger}\widehat{\vek C}_{j} \ebmat,
		\ee
		which yields \eqref{eqn.block-Y-update}.  The proof that these vectors are the minimizers of \eqref{eqn.bl-GMRES-iter-argmin}
		and \eqref{eqn.bl-FOM-iter-argmin} proceeds exactly as in that of Proposition \ref{prop.GF-resid-aug}.  
	\eproof
	
	The behavior of block FOM and GMRES thus can be divided into three cases.
	\begin{enumerate}[{Case} 1]
		\item\label{item.case1} If $\vek H_{j}^{(B)}$ is nonsingular (i.e., the block FOM solution exists) 
			then \eqref{eqn.bl-FOM-iter-argmin} is satisfied exactly, and by augmenting with $L$ 
			columns to expand $\overline{\vek H}_{j-1}^{(B)}$, to 
			${\vek H}_{j}^{(B)}$, the $(j-1)$st GMRES least squares problem becomes a nonsingular linear system.
		\item\label{item.case2} If $\vek H_{j}^{(B)}$ is singular with rank $(j-1)L + r$ with $1\leq r < L$, then 
		the linear system produced by the
			augmentation of $\overline{\vek H}_{j-1}^{(B)}$ produces a better minimizer than $\vek X_{j-1}^{(G)}$ from 
			\eqref{eqn.bl-FOM-iter-argmin}, but it is not exactly solvable.
			This corresponds to only an 
			$r$-dimensional subspace of $\CR(\vek V_{j})$ contributing to the block GMRES minimization at iteration $j$.
		\item\label{item.case3} If $\vek H_{j}^{(B)}$ is singular with rank $(j-1)L$, 
			then the situation is analogous to that described in 
			Theorem \ref{thm.brown-hj-singular}.  We have $\widetilde{\vek X}_{j}^{(F)} = \vek X_{j}^{G} = \vek X_{j-1}^{(G)}$,
			and augmentation of $\overline{\vek H}_{j-1}^{(B)}$ produces no improvement.
	\end{enumerate}
	We note that Case \ref{item.case2} is unique to the block setting and represents a block generalization of the concept of GMRES stagnation,  
	where only an $r$-dimensional subspace of $\CR(\vek V_{j})$ (with $r < \rank \vek V_{j} = L$) contributes to the 
	minimization of the residual at step $j$.  
	We direct the reader to the related discussion in \cite{B.1991} about
	ascent directions, though we omit here such an analysis in the interest of manuscript length.
	Before proving these results, we prove some intermediate technical results.
	
	Let us begin by discussing the structure of $\vek Q_{j}^{(j+1)}$.  In this case, as discussed in Lemma \ref{lemma.C}, this matrix has
	a large $(j-1)L\times (j-1)L$ identity matrix in the upper left-hand corner, and a $2L\times 2L$ nontrivial orthogonal transformation
	block in the lower right-hand corner, denoted 
	\be\label{eqn.oth-tran-decomp}
		\widehat{\boldsymbol\CH}_{j} = \bbmat \vek Q_{j}^{(11)} & \vek Q_{j}^{(12)}\\ \vek Q_{j}^{(21)} & \vek Q_{j}^{(22)} \ebmat,
	\ee
	which we note is itself a product of elementary orthogonal transformations, 
	and all four blocks are of size $L\times L$.  
	Because $\widehat{\boldsymbol\CH}_{j}$ is an orthogonal transformation, it
	admits a CS-decomposition (see, e.g., \cite[Theorem 2.5.3]{GV.2013} and more generally for complex matrices \cite{PW.1994} and
	references therein) i.e., there exist unitary matrices 
	$\boldsymbol\CU_{1},\boldsymbol\CU_{2}, \boldsymbol\CV_{1}, \boldsymbol\CV_{2}\in\C^{L\times L}$ and
	diagonal matrices $\boldsymbol\CS,\boldsymbol\CC\in\R^{L\times L}$ with $\boldsymbol\CS~=~\diag\curl{s_{1},\ldots,s_{L}}$ and 
	$\boldsymbol\CC = \diag\curl{c_{1},\ldots, c_{L}}$ such that 
	\bea
			\vek Q_{j}^{(11)}=\boldsymbol\CU_{1}\boldsymbol\CC\boldsymbol\CV_{1},\, \vek Q_{j}^{(12)}=\boldsymbol\CU_{1}\boldsymbol\CS\boldsymbol\CV_{2},\, \vek Q_{j}^{(21)}&=&\boldsymbol\CU_{2}\boldsymbol\CS\boldsymbol\CV_{1},\nn\\ \mand \vek Q_{j}^{(22)}&=&-\boldsymbol\CU_{2}\boldsymbol\CC\boldsymbol\CV_{2}\label{eqn.CS-decomp-Q},
	\eea
	and for $1\leq i\leq L$ we have $s_{i}^{2} + c_{i}^{2} = 1$, i.e., the diagonal entries of $\boldsymbol\CS$ and $\boldsymbol\CC$ are the sines
	and cosines of $L$ angles, $\curl{\theta_{1},\ldots, \theta_{L}}$.  
	We assume that $c_{1} \leq c_{2}\leq \cdots \leq c_{L}$ and it then follows that $s_{1} \geq s_{2}\geq \cdots \geq s_{L}$.	
	Note that in the case of the single-vector Krylov
	methods, $\widehat{\boldsymbol\CH}_{j}\in\C^{2\times 2}$, 
	$\boldsymbol\CU_{1}=\boldsymbol\CU_{2}=\boldsymbol\CV_{1}=\boldsymbol\CV_{2}=1$, and $\boldsymbol\CS=s_{1}$ and $\boldsymbol\CC=c_{1}$
	are the Givens sine and cosine.  Thus this CS-decomposition yields a nice generalization of the Givens sine and cosine
	in the block setting; see, cf. Section \ref{sec.angle-interp} below	.  We can characterize some elements of this CS-decomposition by studying the QR-factorization of 
	$\overline{\vek H}^{(B)}_{j}$ and its relationship to the rank of $\vek H_{j}^{(B)}$.  The proofs that follow often use generalizations
	of elements of proofs in \cite{B.1991}.  
	\blem\label{lemma.Hbar-R-structure}
		Let $\rank \vek H_{j}^{(B)} = (j-1)L + r$ with $1\leq r \leq L$.
Then we can write
\be\label{eqn.Hjb-rep}
	\vek H_{j}^{(B)} = \bbmat \overline{\vek H}_{j-1}^{(B)} & \vek L_{j} \ebmat
\ee
with $\vek L_{j}\in\C^{jL\times L}$
  such that
		\be\label{eqn.Lj-structure}
			\vek L_{j} = \overline{\vek H}_{j-1}^{(B)}\widehat{\vek Y}_{1} + \vek G_{j}\widehat{\vek Y}_{2}
		\ee 
with $\widehat{\vek Y}_{1}\in\C^{(j-1)L\times L}$, $\widehat{\vek Y}_{2}\in\C^{r\times L}$, and
$\vek G_{j}\in\C^{jL\times r}$ having orthonormal columns which are orthogonal to $\CR\prn{\vek H_{j-1}^{(B)}}$.	
		Furthermore, the blocks $\vek Z_{j}$ and $\widehat{\vek H}_{jj}$ from \eqref{eqn.Rj-orthog-trans} have the following representations
		\be\label{eqn.Hbar-R-structure}
			\vek Z_{j} = \vek R_{j-1}\widehat{\vek Y}_{1} \mand \widehat{\vek H}_{jj} = \widehat{\vek M}_{j}\widehat{\vek Y}_{2}
		\ee
		where 
		$\widehat{\vek M}_{j}\in\C^{L\times r}$ so that
		$\widehat{\vek M}_{j}\widehat{\vek Y}_{2}$ is a rank-$r$ outer product. 
	\elem
	\bproof
		We begin as in \cite{B.1991} by observing that the square matrix $\vek H_{j}^{(B)}$ has the form
		\eqref{eqn.Hjb-rep} following from its nested structure and rank.
		Since $\rank \vek H_{j}^{(B)} = (j-1)L + r$, we can represent the columns of $\vek L_{j}$ as linear combinations of vectors 
		coming from $\CR(\overline{\vek H}_{j-1}^{(B)})$ and vectors coming from a subspace of $\CR(\overline{\vek H}_{j}^{(B)})^{\perp}$, from which \eqref{eqn.Lj-structure} follows,
		where $\vek G_{j}\in\C^{jL\times r}$ has orthonormal columns such that 
		$\CR(\overline{\vek H}_{j-1}^{(B)}) \perp \CR(\vek G_{j})$ and
		$\CR(\vek H_{j}^{(B)})=\CR\prn{\bbmat \overline{\vek H}_{j-1}^{(B)} & \vek G_{j} \ebmat}$. Thus we can write
		\be\nn
			\overline{\vek H}_{j}^{(B)} = \bbmat \overline{\vek H}_{j}^{(B)} & \vek L_{j} \\ & \vek H_{j+1,j} \ebmat = \bbmat \overline{\vek Q}_{j-1}\bbmat \vek R_{j-1}\\\vek 0 \ebmat & \overline{\vek Q}_{j-1}\bbmat \vek R_{j-1}\\\vek 0 \ebmat\overline{\vek Y}_{1} + \vek G_{j}\overline{\vek Y}_{2} \\ & \vek H_{j+1,j} \ebmat ,
		\ee
		and we have that 
		\bea
			\overline{\vek R}_{j} &=& \vek Q_{j}^{(j+1)}\overline{\vek Q}_{j-1}^{\,\ast}\bbmat \overline{\vek Q}_{j-1}\bbmat \vek R_{j-1}\\\vek 0 \ebmat & \overline{\vek Q}_{j-1}\bbmat \vek R_{j-1}\\\vek 0 \ebmat\widehat{\vek Y}_{1} + \vek G_{j}\widehat{\vek Y}_{2} \\ & \vek H_{j+1,j} \ebmat\nn\\
			&=& \vek Q_{j}^{(j+1)}\bbmat \bbmat \vek R_{j-1}\\\vek 0 \ebmat & \bbmat \vek R_{j-1}\\\vek 0 \ebmat\widehat{\vek Y}_{1} + \overline{\vek Q}_{j-1}^{\,\ast}\vek G_{j}\widehat{\vek Y}_{2} \\ & \vek H_{j+1,j} \ebmat\nn.
		\eea
		Since $\overline{\vek Q}_{j-1}\in\C^{(j+1)L\times (j+1)L}$, its columns form an orthonormal basis for $\C^{(j+1)L}$.
		However,  from the upper triangular structure of $\overline{\vek R}_{j-1}$, we know we can partition the columns 
		of $\overline{\vek Q}_{j-1}\in\C^{(j+1)L\times (j+1)L}$ such that the first $(j-1)L$ columns form a basis of 
		$\CR(\overline{\vek H}_{j-1}^{(B)})$ and the remaining columns form a basis for 
		$\CR(\overline{\vek H}_{j-1}^{(B)})^{\perp}$, of which $\CR(\vek G_{j})$ is a subspace.  Thus we can write
		\be\nn
			\overline{\vek Q}_{j-1}^{\,\ast}\vek G_{j} = \bbmat \vek 0\\ \widehat{\vek M}_{j} \ebmat
		\ee
		with $\widehat{\vek M}_{j}\in\C^{L\times r}$ which yields
		\bea
			\overline{\vek R}_{j} &=& \vek Q_{j}^{(j+1)}\bbmat \bbmat \vek R_{j-1}\\\vek 0 \ebmat & \bbmat \vek R_{j-1}\\\vek 0 \ebmat\widehat{\vek Y}_{1} + \bbmat \vek 0\\ \widehat{\vek M}_{j} \ebmat\widehat{\vek Y}_{2} \\ & \vek H_{j+1,j} \ebmat\nn\\
			&=&\vek Q_{j}^{(j+1)}\bbmat \bbmat \vek R_{j-1}\\\vek 0 \ebmat & \bbmat \vek R_{j-1}\widehat{\vek Y}_{1}\\\widehat{\vek M}_{j}\widehat{\vek Y}_{2} \ebmat  \\ & \vek H_{j+1,j} \ebmat\nn
		\eea
		After some simplifications, both the identities for $\vek Z_{j}$ and $\widehat{\vek H}_{jj}$ have been proven.
	\eproof
	\bcor\label{cor.M-is-orthog-Y-is-Tri}
		The representations in \eqref{eqn.Hbar-R-structure} are not unique, and there always exists one such representation
		such that $\widehat{\vek M}_{j}$ has orthonormal columns and $\widehat{\vek Y}_{2}$ is upper triangular.
	\ecor
	\bproof
		Let $\widehat{\vek Y}_{2} = \vek Q_{\widehat{\vek Y}_{2}}\vek R_{\widehat{\vek Y}_{2}}$ be the QR-factorization.
		With the updates $\vek G_{j}\leftarrow \vek G_{j}\vek Q_{\widehat{\vek Y}_{2}}$ and 
		$\widehat{\vek Y}_{2}\leftarrow \vek R_{\widehat{\vek Y}_{2}}$,  \eqref{eqn.Lj-structure} still holds
		with $\vek G_{j}$ still having orthonormal columns.
		With the updates $\widehat{\vek M}_{j} \leftarrow \vek Q_{\widehat{\vek M}_{j}}$ and 
		$\widehat{\vek Y}_{2}\leftarrow \vek R_{\widehat{\vek M}_{j}}\widehat{\vek Y}_{2}$,  
		\eqref{eqn.Hbar-R-structure} still holds.  Thus we have have demonstrated the non-uniqueness of 
		\eqref{eqn.Hbar-R-structure} and that $\widehat{\vek M}_{j}$ and $\widehat{\vek Y}_{2}$ with the
		structure we sought always exist.  
	\eproof
	Henceforth, we assume that $\widehat{\vek M}_{j}$ has orthonormal columns and that $\widehat{\vek Y}_{2}$
	is upper triangular.
	Lemma \ref{lemma.Hbar-R-structure} and Corollary \ref{cor.M-is-orthog-Y-is-Tri} illuminates 
	various properties of the CS-decomposition of $\widehat{\boldsymbol\CH}_{j}$.  We note here that for any $1\leq m \leq L$
	and a matrix $\boldsymbol \CA\in\C^{L\times m}$ with orthonormal columns, that 
	$\boldsymbol \CA^{\perp}\in\C^{L\times (L-r)}$ (a notation we abuse) is some matrix which has orthonormal columns 
	spanning $\CR(\boldsymbol\CA)^{\perp}$ \emph{whose exact structure is determined by the context in which it is used.}
	Furthermore, let $\boldsymbol \CU(\cdot)$ refers to the $\CU$-factor of the singular value decomposition of the argument.
	
	\blem\label{lemma.Horthtrans-props}
		The orthogonal transformation $\widehat{\boldsymbol\CH}_{j}$ with CS-decomposition described in 
		\eqref{eqn.CS-decomp-Q} has the following properties,
		\begin{enumerate}[\rm (I)]
			\item\label{item.Q12} $\vek Q_{j}^{(12)} = \vek N_{j}^{-\ast}\vek H_{j+1,j}^{\ast}$, and it is lower triangular.
			\item\label{item.U1-structure} $\boldsymbol \CU_{1} = \boldsymbol \CU(\vek N_{j}^{-\ast}\vek H_{j+1,j}^{\ast})$.
			\item\label{item.Q21-rank} $\rank \vek Q_{j}^{(12)} = \rank\vek Q_{j}^{(21)} = L$, i.e., they are nonsingular.
			\item\label{item.Q11-rank} $\rank \vek Q_{j}^{(11)} = \rank \vek Q_{j}^{(22)} = r$.
			\item\label{item.V1-structure} $\boldsymbol \CV_{1} = \bbmat \widehat{\vek M}_{j}\boldsymbol\CQ & \widehat{\vek M}_{j}^{\perp} \ebmat$ where $\boldsymbol\CQ\in\C^{r\times r}$ is unitary.
		\end{enumerate}
	\elem
	\bproof
		Observing that
		\be\label{eqn.action-of-orth-trans}
			\bbmat \vek Q_{j}^{(11)} & \vek Q_{j}^{(12)} \\ \vek Q_{j}^{(21)} & \vek Q_{j}^{(22)} \ebmat \bbmat \widehat{\vek H}_{jj}\\ \vek H_{j+1,j} \ebmat = \bbmat \vek N_{j} \\ \vek 0 \ebmat \iff \bbmat \prn{\vek Q_{j}^{(11)}}^{\ast} & \prn{\vek Q_{j}^{(21)}}^{\ast} \\ \prn{\vek Q_{j}^{(12)}}^{\ast} & \prn{\vek Q_{j}^{(22)}}^{\ast} \ebmat \bbmat \vek N_{j} \\ \vek 0 \ebmat = \bbmat \widehat{\vek H}_{jj}\\ \vek H_{j+1,j} \ebmat
		\ee
		and that the right-hand equation of \eqref{eqn.action-of-orth-trans} yields
		\be\label{eqn.transpose-Q-action}
			\prn{\vek Q_{j}^{(11)}}^{\ast}\vek N_{j} = \widehat{\vek H}_{jj} = \widehat{\vek M}_{j}\widehat{\vek Y}_{2},\mand \prn{\vek Q_{j}^{(21)}}^{\ast}\vek N_{j} = \vek H_{j+1,j}.
		\ee
		Since we assume no breakdown of the block Arnoldi method, we know that $\vek N_{j}$ is nonsingular and we can see
		that $\prn{\vek Q_{j}^{(21)}}^{\ast} = \vek H_{j+1,j}\vek N_{j}^{-1}$ which yields Property \ref{item.Q12}.  
		This automatically proves Property \ref{item.U1-structure} as well.		
		This also implies
		that $\vek Q_{j}^{(12)}$ is nonsingular (i.e., rank $L$).  From \eqref{eqn.CS-decomp-Q}, we know $\vek Q_{j}^{(21)}$
		and $\vek Q_{j}^{(12)}$ have the same singular values which completes 
		the proof of Property \ref{item.Q21-rank}.  The first equation
		in \eqref{eqn.transpose-Q-action} can be transformed to $\prn{\vek Q_{j}^{(11)}}^{\ast} = \widehat{\vek M}_{j}\widehat{\vek Y}_{2}\vek N_{j}^{-1}$ implying that $\CR\prn{\prn{\vek Q_{j}^{(11)}}^{\ast}}\subseteq\CR(\widehat{\vek M}_{j})$.  
		We know that $\widehat{\vek Y}_{2}$
		is full rank from how it was constructed, thus $\CR\prn{\prn{\vek Q_{j}^{(11)}}^{\ast}}=\CR\prn{\widehat{\vek M}_{j}}$.
		This yields Property \ref{item.Q11-rank}, since from \eqref{eqn.CS-decomp-Q} we know that $\vek Q_{j}^{(11)}$ and 
		$\vek Q_{j}^{(22)}$ also share the same singular values.  From \eqref{eqn.CS-decomp-Q}, we know that 
		$\prn{\vek Q_{j}^{(11)}}^{\ast} = \boldsymbol\CV_{1}\boldsymbol\CS\boldsymbol\CU_{1}^{\ast}$.
		This implies Property \ref{item.V1-structure} due to the assumed ordering of the singular values contained
		in $\boldsymbol\CS$.  
	\eproof
	Lemma \ref{lemma.Hbar-R-structure} also allows us to describe the structure of the orthogonal transformation
	$\widehat{\vek Q}_{j}^{(b)}$, the non-trivial block of $\widehat{\vek Q}_{j}^{(j)}$.
	\blem\label{lemma.Qhat-structure}
		We have that
		\be\label{eqn.Qjb}
			\widehat{\vek Q}_{j}^{(b)} = \bbmat \widehat{\vek M}_{j} & \widehat{\vek M}_{j}^{\perp} \ebmat^{\ast},
		\ee
		so that we then can write
		\be\label{eqn.Njhat-structure}
			\widehat{\vek N}_{j} = \bbmat \widehat{\vek Y}_{2} \\ \vek 0_{(L-r)\times L} \ebmat.
		\ee
	  \elem
	  \bproof
	  	This follows directly from the assumptions on $\widehat{\vek M}_{j}$ (orthogonal columns) and $\widehat{\vek Y}_{2}$ (upper triangular).
	  \eproof
	  \bcor\label{cor.Cj-rank}
	  	It follows directly that $\rank \vek C_{j} = \rank\widehat{\vek N}_{j}$.
	  \ecor
	  \bproof
	  	The combination of Lemma \ref{lemma.C} with Property \ref{item.Q11-rank} of Lemma \ref{lemma.Horthtrans-props} yields the result.
	  \eproof
	  We have now collected sufficient intermediate results to develop our main results. As in the single-vector Krylov method case, 
	  the rank of $\vek H_{j}^{(B)}$ is intimately related with the solution of the block GMRES least-squares problem \eqref{eqn.block-GMRES-subprob}.
	  The following theorem is a generalization of \cite[Theorem 3.1]{B.1991}, although we frame it a bit differently in this case.
	\bthm\label{thm.H-sing}
		 The matrix $\vek H_{j}^{(B)}$ is singular with 
		$\rank\vek H_{j}^{(B)}=(j-1)L + r$ with $r<L$ if and only if the $j$th block GMRES update $\vek S_{j}^{(G)}$ 
		is such that
		\be\nn
			\dim\prn{\CR\prn{\vek S_{j}^{(G)}}\cap\CR\prn{\vek V_{j}}} = r.
		\ee 
	\ethm	
	\bproof
		Let us first assume that $\vek H_{j}^{(B)}$ is singular with rank $(j-1)L + r$.  It follows then from 
		Lemmas \ref{lemma.block-GF-resid-aug} that
		\be\label{eqn.blGMRES-update-splitting}
			\vek S_{j}^{(G)} = \vek W_{j}\bbmat\
			\vek R_{j-1}^{-1}\vek Z_{j}\vek N_{j}^{-1}\vek C_{j}\\ \vek N_{j}^{-1}\vek C_{j}\ebmat = \vek W_{j-1}\prn{\vek R_{j-1}^{-1}\vek Z_{j}\vek N_{j}^{-1}\vek C_{j}} + \vek V_{j}\prn{\vek N_{j}^{-1}\vek C_{j}}.
		\ee
		From Corollary \ref{cor.Cj-rank} it follows that the rank of $\vek C_{j}$ (and thus also $\vek N_{j}^{-1}\vek C_{j}$) is
		$r$.  Let
		\be\label{eqn.Pj-newSubspaceBasis}
			\boldsymbol\CP_{j} = \bbmat \vek p_{1} & \vek p_{2} & \cdots & \vek p_{r} \ebmat\in\C^{L\times r}
		\ee		
		be the
		matrix with orthonormal columns spanning $\CR(\vek N_{j}^{-1}\vek C_{j})$.
		It follows directly then that the vectors in $\CR\prn{\vek S_{j}^{(G)}}$ only have non-trivial intersection with an
		$r$-dimensional subspace of $\CR\prn{\vek V_{j}}$, namely the subspace $\CR(\vek V_{j}\boldsymbol\CP_{j})$.
		
		Now assume that at the $j$th iteration of block GMRES, the span of the columns of the update $\vek S_{j}^{(G)}$
		has an $r$-dimensional non-trivial intersection with $\CR(\vek V_{j})$.  This implies that there exists
		$\boldsymbol\CP_{j}$ of the form \eqref{eqn.Pj-newSubspaceBasis} such that 
		$\CR\prn{\vek S_{j}^{(G)}}\cap \CR\prn{\vek V_{j}} = \CR\prn{\vek V_{j}\boldsymbol\CP_{j}}$.
		It follows again from 
		Lemma \ref{lemma.block-GF-resid-aug} that $\vek S_{j}^{(G)}$ has the form \eqref{eqn.blGMRES-update-decomp}.
		However, this then implies that $\rank \vek N_{j}^{-1}\vek C_{j} = r$.  Since $\vek N_{j}^{-1}$ is invertible,
		it follows that $\rank \vek C_{j} = r$, and from Corollary \ref{cor.Cj-rank} we then have that 
		$\rank\widehat{\vek N}_{j} = r$, and thus $\vek H_{j}^{(B)}$ has rank $(j-1)L + r$.
	\eproof
	We observe here that Theorem \ref{thm.H-sing} and its proof hinge on the structure of $\vek C_{j}$. If $r$ is nonzero, it 
	follows that $\vek C_{j}$ must be nonzero but singular due to Corollary \ref{cor.Cj-rank}.  The only case in which we
	can have total stagnation (i.e., $\vek C_{j} = \vek 0$), then, is when $r=0$.  Thus we state the following corollary,
	which is the block analog of \cite[Theorem 3.1]{B.1991}.
	\bcor\label{cor.Hsing-GMRES-total-stag}
		The matrix $\vek H_{j}^{(B)}$ is singular with $\rank\vek H_{j}^{(B)}=(j-1)L$ 
		if and only if block GMRES has totally stagnated with 
		$\vek X_{j}^{(G)} = \vek X_{j-1}^{(G)}$. 
	\ecor
	
	It follows that if there is a nontrivial $\vek S_{j}^{(G)}$ whose columns come from 
	$\K_{j}(\vek A,\vek F_{0})$ yielding a better minimizer, it can  
	be decomposed into a part coming from $\CR\prn{\vek V_{j}}$ and a corresponding part 
	from $\K_{j-1}(\vek A,\vek F_{0})$ which is completely determined by the correction coming from $\CR\prn{\vek V_{j}}$.
	\blem
		Let $\vek S_{j}^{(\cdot)} = \vek S_{j,1}^{(\cdot)} + \vek S_{j,2}^{(\cdot)}$ where 
		$\vek S_{j,1}^{(\cdot)}(:,i)\in \K_{j-1}(\vek A,\vek F_{0})$, 
		$\vek S_{j,2}^{(\cdot)}(:,i)\in\CR\prn{\vek V_{j}}$ for $1\leq i\leq L$, and 
		$(\cdot)$ stands for either $(G)$ or $(F)$.  
		Then $\vek S_{j,1}^{(\cdot)} = \boldsymbol{\mathfrak{N}}_{j}^{(\cdot)}\vek S_{j,2}^{(\cdot)}$
		where $\boldsymbol{\mathfrak{N}}_{j}^{(\cdot)}$ is a nilpotent operator such that 
		\be\nn
			\boldsymbol{\mathfrak{N}}_{j}^{(\cdot)}: \CR\prn{\vek V_{j}} \rightarrow \K_{j-1}(\vek A,\vek F_{0})\mand \boldsymbol{\mathfrak{N}}_{j}^{(\cdot)}: \CR\prn{\vek V_{j}}^{\perp} \rightarrow \curl{\vek 0},
		\ee
	\elem
	i.e., $\CR\prn{\boldsymbol{\mathfrak{N}}_{j}^{(\cdot)}} = \K_{j-1}(\vek A,\vek F_{0})$, 
	and $\calN\prn{\boldsymbol{\mathfrak{N}}_{j}^{(\cdot)}} = \CR\prn{\vek V_{j}}^{\perp}$.
	\bproof
		We prove only for the case $(\cdot)=(G)$, as both proofs proceed in the same way.  
		From \eqref{eqn.blGMRES-update-splitting}, we see that 
		\bea
			\vek S_{j,1}^{(G)} &=& \vek W_{j-1}\vek R_{j-1}^{-1}\vek Z_{j}\vek N_{j}^{-1}\vek C_{j}\nn\\ 
							&=& \vek W_{j-1}\vek R_{j-1}^{-1}\vek Z_{j}\vek V_{j}^{\ast}\vek V_{j}\vek N_{j}^{-1}\vek C_{j}\nn\\ 
							&=& \vek W_{j-1}\vek R_{j-1}^{-1}\vek Z_{j}\vek V_{j}^{\ast}\vek S_{j,2}^{(G)}\nn\\ 
							&=& \vek W_{j-1}\widehat{\vek Y}_{j}\vek V_{j}^{\ast}\vek S_{j,2}^{(G)}\nn.
		\eea
		Assigning $\boldsymbol{\mathfrak{N}}_{j}^{(G)} = \vek W_{j-1}\widehat{\vek Y}_{j}\vek V_{j}^{\ast}$, one can easily
		check that it satisfies the statements of the lemma.
	\eproof
	The following theorem is a generalization of \cite[Theorem 3.3]{B.1991}.
	\bthm\label{thm.r-dim-update-bl-GMRES-FOM}
		The span of the columns of $\widetilde{\vek S}_{j}^{(F)}$ has a non-trivial intersection with exactly an $r$-dimensional 
		subspace of $\CR(\vek V_{j})$ if and only if the same is true of ${\vek S}_{j}^{(G)}$. 
	\ethm
	\bproof
		We begin with the assumption
		that ${\vek S}_{j}^{(G)}$ has this property.  We know from Theorem \ref{thm.H-sing} that this implies 
		$\rank \vek H_{j}^{(G)} = (j+1)L + r$ and that $\rank \widehat{\vek N}_{j} = r$.  It follows then that
		$\widehat{\vek N}_{j}^{\dagger}$ has a dimension $L-r$ null space.\footnote[4]{Because we know that $\widehat{\vek N}_{j}$
		is upper triangular with an $(L-r)\times (L-r)$ zero block in the bottom right-hand corner, it follows that 
		$\calN(\widehat{\vek N}_{j}^{\dagger}) = \CR(\widehat{\vek N}_{j})^{\perp} = 
		\Span\curl{\vek e_{r+1}^{[L]}, \vek e_{r+2}^{[L]},
		\ldots, \vek e_{L}^{[L]}}$.  Thus we can write 
		$\widehat{\vek N}_{j}^{\dagger} = \bbmat \boldsymbol\ast_{L\times r} & \vek 0_{L\times (L-r)} \ebmat$.}
		From Lemma \ref{lemma.block-GF-resid-aug}, we can write 
		\be\nn
			\widetilde{\vek S}_{j}^{(F)} = \vek W_{j-1}\prn{-\vek R_{j-1}^{-1}\vek Z_{j}\widehat{\vek N}_{j}^{\dagger}\widehat{\vek C}_{j}} + \vek V_{j}\prn{\widehat{\vek N}_{j}^{\dagger}\widehat{\vek C}_{j}}.
		\ee
		Since we know that $\widehat{\vek C}_{j}$ is nonsingular, it follows that 
		$\rank \widehat{\vek N}_{j}^{\dagger}\widehat{\vek C}_{j} = r$.  Thus, using the same argument used at the
		end of the proof of Theorem \ref{thm.H-sing} it follows that $\widetilde{\vek S}_{j}^{(F)}$ 
		only has a non-trivial intersection with in an $r$-dimensional subspace of $\CR(\vek V_{j})$.
		
		For the other direction, we simply carry out the same steps but in reverse order.
	\eproof
	\bcor\label{cor.bl-GMRES-stag-FOM-stag}
		Block GMRES at iteration $j$ totally stagnates with $\vek X_{j}^{(G)} = \vek X_{j-1}^{(G)}$ if and only if 
		$\widetilde{\vek X}_{j}^{(F)} = \vek X_{j-1}^{(G)}$.
	\ecor
	\bproof
		This corresponds to the case $r=0$ for Theorem \ref{thm.r-dim-update-bl-GMRES-FOM}.
	\eproof
	
	We now show that the case of partial stagnation of block GMRES
	(as defined at the beginning of Section \ref{section.main-results}) is actually just a special
	case of Theorem \ref{thm.H-sing}, and is not really of special interest with respect to this analysis
	\bthm\label{thm.partial-stag}
		Block GMRES suffers a partial stagnation at iteration $j$ of the form \eqref{eqn.partial-stag} if and only if 
		$0 < \rank \vek C_{j} \leq r$ where $r=|\overline{\I}|$ such that for all 
		$i\in\I$ the $i$th column of $\vek C_{j}$ is the zero vector.
	\ethm
	\bproof
		Let us first assume that the columns of $\vek C_{j}$ corresponding to indices in $\I$ are zero but that 
		$\vek C_{j}\neq\vek 0$.  Then $\rank \vek C_{j} \leq L-|\I|$. Furthermore, since 
		\be\nn
			\vek S_{j}^{(G)} = \vek W_{j}\bbmat\vek R_{j-1}^{-1}\vek Z_{j}\\ \vek I\ebmat\vek N_{j}^{-1}\vek C_{j},
		\ee
		for $i\in\I$, if we look at the $i$th column of $\vek S_{j}^{(G)}$, we see that 
		\be\label{eqn.blGMRES-update-decomp}
			\vek S_{j}^{(G)}\vek e_{i}^{[L]} = \vek W_{j}\bbmat\vek R_{j-1}^{-1}\vek Z_{j}\\ \vek I\ebmat\vek N_{j}^{-1}\vek C_{j}\vek e_{i}^{[L]} = \vek W_{j}\bbmat\vek R_{j-1}^{-1}\vek Z_{j}\\ \vek I\ebmat\vek N_{j}^{-1}\vek 0 = \vek 0.
		\ee
		The first direction is thus proven.
		
		Now we assume that partial stagnation occurs at the $j$th iteration where for each $i\in\I$, 
		$\vek X_{j-1}^{(G)}\vek e_{i}^{[L]} = \vek X_{j}^{(G)}\vek e_{i}^{[L]}$.  
		This implies that $\vek S_{j}^{(G)}\vek e_{i}^{[L]}=\vek 0$
		for all $i\in\I$.  Specifically, this implies that $\vek V_{j}\vek N_{j}^{-1}\vek C_{j}\vek e_{i}^{[L]} = \vek 0$.
		Because we assume that $\vek V_{j}$ is full rank and $\vek N_{j}^{-1}$ is nonsingular, it follows that 
		$\vek C_{j}\vek e_{i}^{[L]} = \vek 0$, which proves the other direction.
	\eproof
	Now we also state the block analog to Proposition \ref{prop.gmres-single-stag}.
	\bthm\label{thm.bl-GMRES-FOM-trig-rel}
		Suppose that $\vek H_{j}^{(B)}$ is nonsingular.  
		Then at iteration $j$ we have the following relationship between the approximations 
		produced by block GMRES and block FOM,
		\be\label{eqn.bl-GMRES-FOM-trig-rel}
			\vek X_{j}^{(G)} = \vek X_{j}^{(F)}\prn{\widehat{\vek C}_{j}^{-1}\boldsymbol\CQ\boldsymbol\CC^{2}\boldsymbol\CQ^{\ast}\widehat{\vek C}_{j}} + \vek X_{j-1}^{(G)}\prn{\widehat{\vek C}_{j}^{-1}\boldsymbol\CQ\boldsymbol\CS^{2}\boldsymbol\CQ^{\ast}\widehat{\vek C}_{j}}.
		\ee
	\ethm
	\bproof
		Since $\vek H_{j}^{(B)}$ is nonsingular, we have that $\widehat{\vek N}_{j}$ and $\vek C_{j}$ are nonsingular,
		and the block FOM approximation $\vek X_{j}^{(F)}$ (and thus also $\vek Y_{j}^{(F)}$) exists. 
		From the proof of Lemma \ref{lemma.block-GF-resid-aug}, we have then that
		\bea
			\prn{\vek Y_{j}^{(F)} - \bbmat\vek Y_{j-1}^{(G)}\\ \vek 0\ebmat} = \bbmat\vek R_{j-1}^{-1}\vek Z_{j}\\ \vek I \ebmat\widehat{\vek N}_{j}^{-1}\widehat{\vek C}_{j}&\mand & \nn\\ \prn{\vek Y_{j}^{(G)} - \bbmat\vek Y_{j-1}^{(G)}\\ \vek 0\ebmat} = \bbmat\vek R_{j-1}^{-1}\vek Z_{j}\\ \vek I \ebmat{\vek N}_{j}^{-1}{\vek C}_{j}&&\nn.
		\eea
		Because in this case, everything is invertible, we see that 
		\be\label{eqn.block-Y-relation}
			\prn{\vek Y_{j}^{(F)} - \bbmat\vek Y_{j-1}^{(G)}\\ \vek 0\ebmat}\widehat{\vek C}_{j}^{-1}\widehat{\vek N}_{j}\vek N_{j}^{-1}\vek C_{j} = \vek Y_{j}^{(G)} - \bbmat\vek Y_{j-1}^{(G)}\\ \vek 0\ebmat.
		\ee
		We can now simplify $\widehat{\vek C}_{j}^{-1}\widehat{\vek N}_{j}\vek N_{j}^{-1}\vek C_{j}$ using 
		Lemmas \ref{lemma.Chat}, \ref{lemma.C}, and \ref{lemma.Horthtrans-props}.  We note that in the case that the
		block FOM approximation exists, we have that $\widehat{\vek N}_{j} = \widehat{\vek Y}_{j}$ is upper triangular
		and nonsingular, $\boldsymbol\CV_{1} = \widehat{\vek M}_{j}\boldsymbol\CQ$, and $\vek Q_{j}^{(b)} = \widehat{\vek M}_{j}^{\ast}$.  
		It follows then that
		\bea
			\widehat{\vek C}_{j}^{-1}\widehat{\vek N}_{j}\vek N_{j}^{-1}\vek C_{j} & = & \widehat{\vek C}_{j}^{-1}\widehat{\vek Y}_{j} \prn{\prn{\vek Q_{j}^{(11)}}^{-\ast}\widehat{\vek M}_{j}\widehat{\vek Y}_{j} }^{-1}\vek Q_{j}^{(11)}\widetilde{\vek C}_{j}\nn\\
			 & = & \widehat{\vek C}_{j}^{-1}\widehat{\vek Y}_{j} \prn{\widehat{\vek Y}_{j}^{-1}\widehat{\vek M}_{j}^{\ast}\prn{\vek Q_{j}^{(11)}}^{\ast} }\vek Q_{j}^{(11)}\widetilde{\vek C}_{j}\nn\\
			 & = & \widehat{\vek C}_{j}^{-1}\widehat{\vek Y}_{j} {\widehat{\vek Y}_{j}^{-1}\widehat{\vek M}_{j}^{\ast}\prn{\boldsymbol\CU_{1}\boldsymbol\CC\boldsymbol\CV_{1}^{\ast}}^{\ast} }\boldsymbol\CU_{1}\boldsymbol\CC\boldsymbol\CV_{1}^{\ast}\widetilde{\vek C}_{j}\nn\\
			 & = & \widehat{\vek C}_{j}^{-1}\widehat{\vek Y}_{j} {\widehat{\vek Y}_{j}^{-1}\widehat{\vek M}_{j}^{\ast}\prn{\widehat{\vek M}_{j}\boldsymbol\CQ\boldsymbol\CC\boldsymbol\CU_{1}^{\ast}} }\boldsymbol\CU_{1}\boldsymbol\CC\prn{\widehat{\vek M}_{j}\boldsymbol\CQ}^{\ast}\widetilde{\vek C}_{j}\nn\\
			 & = & \widehat{\vek C}_{j}^{-1} {{\boldsymbol\CQ} }\boldsymbol\CC^{2}\prn{\widehat{\vek M}_{j}\boldsymbol\CQ}^{\ast}\widetilde{\vek C}_{j}\nn\\
			 & = & \widehat{\vek C}_{j}^{-1} {{\boldsymbol\CQ} }\boldsymbol\CC^{2}{\boldsymbol\CQ^{\ast}}\widehat{\vek C}_{j}\label{eqn.Csquared-deriv}.
		\eea
		We now insert \eqref{eqn.Csquared-deriv} into \eqref{eqn.block-Y-relation}, multiply both sides by $\vek W_{j}$,
		and perform some algebraic manipulations to get
		\be\nn
			\vek X_{j}^{(G)} = \vek X_{j}^{(F)}\prn{\widehat{\vek C}_{j}^{-1}\boldsymbol\CQ\boldsymbol\CC^{2}\boldsymbol\CQ^{\ast}\widehat{\vek C}_{j}} + \vek X_{j-1}^{(G)}\prn{\vek I - \widehat{\vek C}_{j}^{-1}\boldsymbol\CQ\boldsymbol\CC^{2}\boldsymbol\CQ^{\ast}\widehat{\vek C}_{j}}.
		\ee
		Lastly, we observe that  ${\widehat{\vek C}_{j}^{-1}\boldsymbol\CQ\boldsymbol\CC^{2}\boldsymbol\CQ^{\ast}\widehat{\vek C}_{j}}$ is an 
		eigen-decomposition since 
		$\prn{\boldsymbol\CQ^{\ast}\widehat{\vek C}_{j}}^{-1} = \widehat{\vek C}_{j}^{-1}\boldsymbol\CQ$, and
		we thus can write
		\be\nn
			\vek I - \widehat{\vek C}_{j}^{-1}\boldsymbol\CQ\boldsymbol\CC^{2}\boldsymbol\CQ^{\ast}\widehat{\vek C}_{j} = \widehat{\vek C}_{j}^{-1}\boldsymbol\CQ\prn{\vek I - \boldsymbol\CC^{2}}\boldsymbol\CQ^{\ast}\widehat{\vek C}_{j}.
		\ee
		The result follows by observing that $\vek I - \boldsymbol\CC^{2} = \boldsymbol\CS^{2}$ which follows from \eqref{eqn.CS-decomp-Q}.
	\eproof
	We will return shortly to understand the meaning of the angles associated to these sines and cosines in
	Section \ref{sec.angle-interp}.
	
	\subsection{The case of breakdown in the block Arnoldi process}\label{section.breakdown-case}
	Our discussion of the case of breakdown focuses first upon the behaviors of block GMRES and block FOM at the $j$th iteration 
	in which the block Arnoldi process produces $p$ dependent basis vectors.  For simplicity, we assume that no single system
	has converged but rather that some linear combination of some columns of the solution $\vek X$ is in $\K_{j}(\vek A, \vek F_{0})$.
	Both the block GMRES and block FOM residuals are thus of rank $L-p$.  
	We assume that these $p$ vectors are replaced with $p$ random vectors so that we maintain a block size of $L$.  For the most part, what we have proven thus far holds with 
	little to no alterations, but the reduction of residual rank does have some consequences.
	
	We consider a breakdown at iteration $j$ in which $p$ dependent basis vectors are produced.
	Various strategies have been suggested for replacing dependent vectors in the interest of maintaining
	the block size of $L$.
	The block Arnoldi process produces from $\vek A\vek V_{j}\in\C^{n\times L}$ the block
	\be\nn
		\vek U_{j+1} = \vek A\vek V_{j} - \sum_{i=1}^{j}\vek V_{i}\vek H_{ij}\mwith\vek H_{ij}\in\C^{L\times L},\mand \rank\vek U_{j+1}=L-p.
	\ee
	Then we have the reduced QR-factorization $\vek U_{j+1} = \ddot{\vek V}_{j+1}\ddot{\vek H}_{j+1,j}$ where 
	$\ddot{\vek V}_{j+1}\in\C^{n\times (L-p)}$ has
	columns spanning $\CR(\vek U_{j+1})$, and $\ddot{\vek H}_{j+1,j}\in\C^{(L-p)\times L}$ is upper triangular.  
	To maintain block size, we set $\vek H_{j+1,j} = \bbmat \ddot{\vek H}_{j+1,j} \\ \vek 0\ebmat\in\C^{L\times L}$ 
	and $\vek V_{j+1} = \bbmat \ddot{\vek V}_{j+1} & \boldsymbol\CZ \ebmat \in\C^{n\times L}$ where 
	$\boldsymbol\CZ\in\C^{n\times p}$ are the independent replacement vectors, which have been orthonormalized against all
	of the block Arnoldi vectors.  Thus the columns of $\vek W_{j+1}$ no longer span a Krylov subspace, but they do still
	satisfy the block Arnoldi relation \eqref{eqn.block-arnoldi-relation}.  The iteration continues unabated.  It is observed in,
	e.g., \cite{RS.2006}, that at the iteration in which the breakdown occurs, the least squares problem still has a unique solution.
	From the analysis in this paper, this corresponds to $\vek N_{j}$ still being nonsingular.  
	Furthermore, as we assume that iteration $j$ is the first iteration at which there is a block Arnoldi breakdown, the
	block residual $\vek F_{j-1}^{(G)}$ is full rank; and, thus, so are the $\vek C$-matrices.  	
	Therefore, at iteration $j$, if there has been a block Arnoldi breakdown, 
	all of the results we have proven still hold with no alteration.  The block residual
	$\vek F_{j}^{(G)}$ has rank $L-p$.  
	
	Without loss of generality, let us consider the case that the breakdown at iteration $j$ is the only breakdown.  
	Consider some later iteration $j+k$ with $k>0$.
	As we have replaced all dependent Arnoldi vectors with 
	linearly independent ones, the GMRES least squares problem still has a unique solution.  This implies that $\vek N_{j+k}$
	is still nonsingular.  
	The block GMRES residuals will continue to be rank $L-p$.  Thus, the $\vek C$-matrices
	will be square (as we maintain block size) and rank-deficient.
	However, few of the results rely on the invertibility of 
	these matrices.  Indeed, the only result not valid in this case is Theorem \ref{thm.bl-GMRES-FOM-trig-rel}.  However, we can prove
	a weaker result in this case.
	\bthm\label{thm.weaker-result}
		Suppose at step $j$ there has been a block Arnoldi breakdown with $p$ dependent Arnoldi vectors being generated, 
		and that there are no further breakdowns
		Let these vectors
		be replaced using the procedure described above. Then at iteration $j+k$, if
		$\vek H_{j+k}^{(B)}$ is nonsingular we have that 
		\be\nn
			\vek X_{j}^{(F)} - \vek X_{j}^{(G)} = \vek W_{j}\bbmat\vek R_{j-1}^{-1}\vek Z_{j}\\ \vek I \ebmat\widehat{\vek Y}_{2}^{-1}\boldsymbol\CQ\boldsymbol\CS^{2}\boldsymbol\CQ^{\ast}\widehat{\vek C}_{j}.
		\ee
	\ethm
	\bproof
		We show this by substituting many of the identities we have previously proven, which are still valid in this setting.
		\bea
		\vek Y_{j}^{(F)} - \vek Y_{j}^{(G)} &=& \bbmat\vek R_{j-1}^{-1}\vek Z_{j}\\ \vek I \ebmat\prn{\widehat{\vek N}_{j}^{-1}\widehat{\vek Q}_{j}^{(b)} - \vek N_{j}^{-1}\vek Q_{j}^{(11)}}\widetilde{\vek C}_{j}\nn\\
		&=& \bbmat\vek R_{j-1}^{-1}\vek Z_{j}\\ \vek I \ebmat\prn{\widehat{\vek Y}_{2}^{-1}\widehat{\vek M}_{j}^{\ast} - \widehat{\vek Y}_{2}^{-1}\widehat{\vek M}_{j}^{\ast}\vek Q_{j}^{(11)\ast}\vek Q_{j}^{(11)}}\widetilde{\vek C}_{j}\nn\\
		&=& \bbmat\vek R_{j-1}^{-1}\vek Z_{j}\\ \vek I \ebmat\prn{\widehat{\vek Y}_{2}^{-1}\widehat{\vek M}_{j}^{\ast} - \widehat{\vek Y}_{2}^{-1}\widehat{\vek M}_{j}^{\ast}\boldsymbol\CV_{1}\boldsymbol\CC\boldsymbol\CU_{1}^{\ast}\boldsymbol\CU_{1}\boldsymbol\CC\boldsymbol\CV_{1}^{\ast}}\widetilde{\vek C}_{j}\nn\\
		&=& \bbmat\vek R_{j-1}^{-1}\vek Z_{j}\\ \vek I \ebmat\prn{\widehat{\vek Y}_{2}^{-1}\widehat{\vek M}_{j}^{\ast} - \widehat{\vek Y}_{2}^{-1}\widehat{\vek M}_{j}^{\ast}\widehat{\vek M}_{j}\boldsymbol\CQ\boldsymbol\CC^{2}\boldsymbol\CQ^{\ast}\widehat{\vek M}_{j}^{\ast}}\widetilde{\vek C}_{j}\nn\\
		&=& \bbmat\vek R_{j-1}^{-1}\vek Z_{j}\\ \vek I \ebmat\widehat{\vek Y}_{2}^{-1}\prn{\vek I - \boldsymbol\CQ\boldsymbol\CC^{2}\boldsymbol\CQ^{\ast}}\widehat{\vek C}_{j}.\nn
		\eea
		On then performs a bit of algebra and multiplies both sides by $\vek W_{j}$ to get the result.
	\eproof
Although this result is less satisfying that Theorem \ref{thm.bl-GMRES-FOM-trig-rel}, as it does not generalize 
Proposition \ref{prop.gmres-single-stag}, it still yields valuable information about the relationship of the block FOM and block
GMRES iterates in the case that breakdown has occurred.  We see that if the angles represented by the sines contained on the
diagonal of $\boldsymbol\CS$ are small, this implies that block FOM and block GMRES in this scenario produce iterations
which are not far from one another.  We must thus now clarify the precise significance of these angles to complete our
analysis.
\subsection{Principal angles between the range of $\vek F_{j-1}^{(G)}$ and $\vek A\K_{j}(\vek A,\vek F_{0})$}\label{sec.angle-interp}
In this section, we show that the angles represented by the sines and cosines from the CS-decomposition 
of \eqref{eqn.oth-tran-decomp} which appear in \eqref{eqn.bl-GMRES-FOM-trig-rel} are the principal angles between the previous
residual and the current residual constraint space.  

In \cite{Eiermann2001}, many geometric properties of single-vector projection methods were analyzed.  In particular,
the authors discussed minimum residual projection methods such as GMRES.  In that paper, the authors show that the angle
represented by the Givens sine and cosine calculated at iteration $j$ of GMRES is actually the principal angle between
the $(j-1)$st GMRES residual and the $j$th constraint space.  In essence, the closeness of this angle to zero indicates how much of
the $(j-1)$st residual lies in the $j$th constraint space, and will thus be eliminated by the projection at iteration $j$.  If the angle is
near $\frac{\pi}{2}$, however, then the Givens cosine $c_{j}$ is close to $0$ and we have near stagnation, since almost
none of the $j-1$st residual lies in the new constraint space and thus there will not be much improvement from the 
projection at iteration $j$.

To illuminate the meaning of these angles in the block setting, we generalize some results from the single-vector GMRES
case.  Following \cite{Eiermann2001}, we represent $\vek A\K_{j}(\vek A,\vek F_{0})$ with a specific, useful basis.
The columns of $\vek W_{j+1}$ form an orthonormal basis for $\K_{j+1}(\vek A,\vek F_{0})$, and it follows from the 
block Arnoldi relation \eqref{eqn.block-arnoldi-relation} that the columns of $\vek W_{j+1}\overline{\vek H}_{j}^{(B)}$ form
a non-orthonormal basis for  $\vek A\K_{j}(\vek A,\vek F_{0})$, and using the QR-factorization 
$\overline{\vek H}_{j}^{(B)} = \overline{\vek Q}_{j}^{\ast}\overline{\vek R}_{j}$, we see that the columns of
$\vek W_{j+1}\overline{\vek Q}_{j}^{\ast}$ forms another orthonormal basis of $\vek A\K_{j}(\vek A,\vek F_{0})$. From the equation
\be\nn
	\vek F_{0} = \vek W_{j+1}\vek E_{1}^{[j+1]}\vek S_{0} = \vek F_{0} = \vek W_{j+1}\overline{\vek Q}_{j}^{\ast}\overline{\vek Q}_{j}\vek E_{1}^{[j+1]}\vek S_{0},
\ee
we see that $\overline{\vek Q}_{j}\vek E_{1}^{[j+1]}\vek S_{0}$ is a representation of $\vek F_{0}$ in that basis.  This leads to a 
generalization of, e.g., \cite[Equation 6.48]{Saad.Iter.Meth.Sparse.2003}, that the Givens sines can be used to cheaply update
the GMRES residual norm.  We note that following from the block partitioning of the orthogonal transformation in 
\eqref{eqn.oth-tran-decomp}, we can write
\be\label{eqn.gen-orth-trans-struct}
	\vek Q_{i}^{(j+1)} = \bbmat \vek I_{(i-1)L} && \\ &\begin{matrix} \vek Q_{i}^{(11)} & \vek Q_{i}^{(12)}\\ \vek Q_{i}^{(21)} & \vek Q_{i}^{(22)} \end{matrix}& \\ && \vek I_{(j-i)L} \ebmat.
\ee
Then we have the following.
\blem
	The representation $\overline{\vek Q}_{j}\vek E_{1}^{[j+1]}\vek S_{0}$  of $\vek F_{0}$ has the following structure,
	\be\label{eqn.resid-vec-struct}
		\overline{\vek Q}_{j}\vek E_{1}^{[j+1]}\vek S_{0} = \bbmat \vek Q_{1}^{(11)} \\ \vek Q_{2}^{(11)}\vek Q_{1}^{(21)} \\ \vek Q_{3}^{(11)}\vek Q_{2}^{(21)}\vek Q_{1}^{(21)} \\ \vdots \\ \vek Q_{j-1}^{(11)} \prod_{i=1}^{j-2}\vek Q_{i}^{(21)} \\  \vek Q_{j}^{(11)}\prod_{i=1}^{j-1}\vek Q_{i}^{(21)} \\ \prod_{i=1}^{j}\vek Q_{i}^{(21)}  \ebmat\vek S_{0}.
	\ee
\elem
\bproof
	This follows from the fact that $\overline{\vek Q}_{j} = \prod_{i=1}^{j}\vek Q_{j-i+1}^{(j+1)}$  and the structure of the
	orthogonal transformations in \eqref{eqn.gen-orth-trans-struct}.
\eproof
Let us denote with $\measuredangle\prn{\CU_{1},\CU_{2}}$ the set of principle angles between the subspaces 
$\CU_{1}$ and $\CU_{2}$.
Following \cite{Eiermann2001}, we can compute a product of matrices whose singular values are the sines of the principal 
angles $\measuredangle\prn{\CR\prn{\vek F_{0}},\vek A\K_{j}(\vek A,\vek F_{0})}$.
\blem
	The principal angles $\measuredangle\prn{\CR\prn{\vek F_{0}},\vek A\K_{j}(\vek A,\vek F_{0})}$ are the singular values of the
	product $\prod_{i=1}^{j}\vek Q_{j-i+1}^{(21)}$.
\elem
\bproof
	We have the equalities
	\bea
		\measuredangle\prn{\CR\prn{\vek F_{0}},\vek A\K_{j}(\vek A,\vek F_{0})} &=& \measuredangle\prn{\CR(\overline{\vek Q}_{j}\vek E_{1}^{[j+1]}\vek S_{0}),\CR(\vek W_{j+1}\overline{\vek Q}_{j+1}^{\ast}\overline{\vek R}_{j})}\nn\\ 
		&=& \measuredangle\prn{\CR\prn{\overline{\vek Q}_{j}\vek E_{1}^{[j+1]}},\CR\prn{\overline{\vek R}_{j}}}.\nn
	\eea
	Under the assumptions in this paper, $\CR\prn{\overline{\vek R}_{j}} = \CR\prn{\bbmat \vek R_{j}\\ \vek 0 \ebmat}$ is isomorphic 
	with $\C^{(j+1)L}$ (due to the nonsingularity of $\vek R_{j}$) with basis $\curl{\vek e_{1}^{[j+1)L]},\vek e_{2}^{[j+1)L]},\ldots,\vek e_{jL}^{[j+1)L]}}$, i.e., the 
	last $L$ coordinates are zero. It is clear that $\overline{\vek Q}_{j}\vek E_{1}^{[j+1]}$ has orthonormal columns.  Let
	\be\nn
		\overline{\vek Q}_{j}\vek E_{1}^{[j+1]} = \bbmat \boldsymbol{\mathfrak{Q}}_{1}\\ \boldsymbol{\mathfrak{Q}}_{2} \ebmat
	\ee
	be a block partitioning with $\boldsymbol{\mathfrak{Q}}_{1}\in\C^{jL\times L}$ and $\boldsymbol{\mathfrak{Q}}_{2}\in\C^{L\times L}$.
	With this partitioning, $\overline{\vek Q}_{j}\vek E_{1}^{[j+1]}$ admits a skinny CS-decomposition (see, e.g., 
	\cite[Section 2.5.4]{GV.2013}) yielding the simultaneous singular value decompositions
	\be\nn
		\boldsymbol{\mathfrak{C}} = \boldsymbol{\mathfrak{U}}_{1}^{\ast}\boldsymbol{\mathfrak{Q}}_{1}\boldsymbol{\mathfrak{V}}\mand \boldsymbol{\mathfrak{S}} = \boldsymbol{\mathfrak{U}}_{2}^{\ast}\boldsymbol{\mathfrak{Q}}_{2}\boldsymbol{\mathfrak{V}}
	\ee
	with $\boldsymbol{\mathfrak{U}}_{1}\in\C^{jL\times jL}$, $\boldsymbol{\mathfrak{U}}_{2}\in\C^{L\times L}$, $\boldsymbol{\mathfrak{V}}\in\C^{L\times L}$, and 
	\be\nn
		\boldsymbol{\mathfrak{C}} = \bbmat \vek I_{(j-1)L} & \\ & {\rm diag}_{i=1}^{L}\curl{c_{i}} \ebmat\in\C^{jL\times jL}\mand \boldsymbol{\mathfrak{S}} = {\rm diag}_{i=1}^{L}\curl{s_{i}}\in\C^{L\times L}.
	\ee
	Since $\bbmat \vek I_{jL}\\\vek 0_{L} \ebmat$ has orthonormal columns spanning $\CR\prn{\overline{\vek R}_{j}}$, 
	the cosines of the sought-after principal angles are given by the singular values of 
	$\boldsymbol{\mathfrak{Q}}_{1} = \bbmat \vek I_{jL}\\\vek 0_{L} \ebmat^{\ast}\prn{\overline{\vek Q}_{j}\vek E_{1}^{[j+1]}}$,
	i.e., the entries of $\boldsymbol{\mathfrak{C}}$.
	Many of these are trivially one (i.e., $\theta_{i}=0$ for $i=1,2,\ldots (j-1)L$).  However, the $L$ nontrivial angles are also
	represented by their sines in the entries of $\boldsymbol{\mathfrak{S}}$ which are the singular values of 
	$\boldsymbol{\mathfrak{Q}}_{2}$, and this proves the lemma.  
\eproof

Using similar techniques, we can prove the following result.
\bthm
The angles represented by the sines and cosines of the \linebreak CS-decomposition of the 
$j$th orthogonal transformation \eqref{eqn.CS-decomp-Q} are the principal angles between the column space of the
previous block GMRES residual $\vek F_{j-1}^{(G)}$ and the $j$th constraint space $\vek A\K_{j}(\vek A,\vek F_{0})$.
\ethm

\bproof
	As has already been discussed, the columns of
	$\vek W_{i+1}\overline{\vek Q}_{i}^{\ast}\bbmat \vek I_{(i-1)L}\\ \vek 0_{L} \ebmat$ are an orthonormal bases for
	$\vek A\K_{i}(\vek A,\vek F_{0})$ for all $i$.  Let $\boldsymbol\CP_{j-1}$ be the orthogonal projector onto  
	$\vek A\K_{j-1}(\vek A,\vek F_{0})$ which means we can write the 
	$\vek F_{j-1}^{(G)} = (\vek I - \boldsymbol\CP_{j-1})\vek F_{0}$.  Using the orthonormal basis of 
	$\vek A\K_{j-1}(\vek A,\vek F_{0})$, we can write
	\be\nn
		\boldsymbol\CP_{j-1}\vek F_{0} = \vek W_{j}\overline{\vek Q}_{j-1}^{\ast}\bbmat \vek I_{(i-1)L}\\ \vek 0_{L} \ebmat\bbmat \vek I_{(i-1)L}&\vek 0_{L} \ebmat\overline{\vek Q}_{j-1}\vek E_{1}^{[jL]}\vek S_{0}.
	\ee
	It is then straightforward to show that
	\be\nn
		\vek F_{j-1}^{(G)} = (\vek I - \boldsymbol\CP_{j-1})\vek F_{0} = \vek W_{j}\overline{\vek Q}_{j-1}^{\ast}\prn{\vek I - \bbmat \vek I_{(i-1)L}\\ \vek 0_{L} \ebmat\bbmat \vek I_{(i-1)L}&\vek 0_{L} \ebmat}\overline{\vek Q}_{j-1}\vek E_{1}^{[jL]}\vek S_{0}.
	\ee
	Observe that $\prn{\vek I - \bbmat \vek I_{(i-1)L}\\ \vek 0_{L} \ebmat\bbmat \vek I_{(i-1)L}&\vek 0_{L} \ebmat}$ is the
	orthogonal projector onto the last $L$ coordinate directions, i.e., onto 
	${\rm span}\curl{\vek e_{(j-1)L+1}^{[mL]},\vek e_{(j-1)L+2}^{[mL]},\ldots,\vek e_{jL}^{[mL]}}$.  Combining this with
	\eqref{eqn.resid-vec-struct}, we can rewrite 
	\be\nn
		(\vek I - \boldsymbol\CP_{j-1})\vek F_{0} = \vek W_{j}\overline{\vek Q}_{j-1}^{\ast}\bbmat \vek 0_{(j-1)L}\\  \\\prod_{i=1}^{j-1}\vek Q_{i}^{(21)} \ebmat\vek S_{0}
	\ee
	The principal angle calculation can then be simplified,
	\small
	\bea
		\measuredangle\prn{\CR\prn{\vek F_{j-1}^{(G)}}, \vek A\K_{j}(\vek A,\vek F_{0})} & = & \measuredangle\prn{\CR\prn{\vek W_{j}\overline{\vek Q}_{j-1}^{\ast}\bbmat \vek 0_{(j-1)L}\\  \\\prod_{i=1}^{j-1}\vek Q_{i}^{(21)} \ebmat\vek S_{0}},\CR\prn{\vek W_{j+1}\overline{\vek Q}_{j}^{\ast}\overline{\vek R}_{j}}}\nn\\
		& = & \measuredangle\prn{\CR\prn{\vek W_{j+1}\bbmat \vek I_{jL}\\ \vek 0_{L} \ebmat\overline{\vek Q}_{j-1}^{\ast}\bbmat \vek 0_{(j-1)L}\\  \\\prod_{i=1}^{j-1}\vek Q_{i}^{(21)} \ebmat},\CR\prn{\vek W_{j+1}\overline{\vek Q}_{j}^{\ast}\overline{\vek R}_{j}}}\nn\\
		& = & \measuredangle\prn{\CR\prn{\bbmat \overline{\vek Q}_{j-1}^{\ast}\\ \vek 0_{L} \ebmat\bbmat \vek 0_{(j-1)L}\\  \\\prod_{i=1}^{j-1}\vek Q_{i}^{(21)} \ebmat},\CR\prn{\overline{\vek Q}_{j}^{\ast}\bbmat \vek I_{jL}\\ \vek 0_{L} \ebmat}}\nn\\
		& = & \measuredangle\prn{\CR\prn{\bbmat \overline{\vek Q}_{j-1}^{\ast}\bbmat \vek 0_{(j-1)L}\\  \\\prod_{i=1}^{j-1}\vek Q_{i}^{(21)} \ebmat\\ \vek 0_{L} \ebmat},\CR\prn{\overline{\vek Q}_{j}^{\ast}\bbmat \vek I_{jL}\\ \vek 0_{L} \ebmat}}\nn
	\eea
	\normalsize
	Observe now that we can rewrite 
	\footnotesize
	\bea
		\overline{\vek Q}_{j}^{\ast}\bbmat \vek I_{jL}\\ \vek 0_{jL\times L} \ebmat &=& \vek Q_{1}^{(j+1)\ast}\cdots\vek Q_{j}^{(j+1)\ast} \bbmat \vek I_{L} &&& \\ &\vek I_{L}&& \\ &&\ddots& \\ &&&\vek I_{L} \\ \vek 0_{L}  &\vek 0_{L}&\cdots &\vek 0_{L}\ebmat \nn\\
		&=& \vek Q_{1}^{(j+1)\ast}\cdots\vek Q_{j-1}^{(j+1)\ast} \bbmat \vek I_{L} &&&& \\ &\vek I_{L}&& \\ &&\ddots &\\&&&\vek I_{L}& \\ &&&&\vek Q_{j}^{(11)} \\ \vek 0_{L}  &\vek 0_{L}&\cdots&\vek 0_{L} &\vek Q_{j}^{(12)\ast}\ebmat \nn\\
		&=& \bbmat \overline{\vek Q}_{j-1}^{\ast} & \\ & \vek I_{L}  \ebmat\bbmat \vek I_{L} &&&& \\ &\vek I_{L}&& \\ &&\ddots &\\&&&\vek I_{L}& \\ &&&&\vek Q_{j}^{(11)} \\ \vek 0_{L}  &\vek 0_{L}&\cdots&\vek 0_{L} &\vek Q_{j}^{(12)\ast}\ebmat, \nn
	\eea
	\normalsize
	and similarly we have
	\footnotesize
	\be\nn
		\bbmat \overline{\vek Q}_{j-1}^{\ast}\bbmat \vek 0_{(j-1)L}\\  \\\prod_{i=1}^{j-1}\vek Q_{i}^{(21)} \ebmat\\ \vek 0_{L} \ebmat = \bbmat \overline{\vek Q}_{j-1}^{\ast} & \\ & \vek I_{L} \ebmat  \bbmat \vek 0_{(j-1)L} \\ \prod_{i=1}^{j-1}\vek Q_{i}^{(21)} \\ \vek 0_{L} \ebmat
	\ee
	\normalsize
	Thus we have 
	\footnotesize
		\bea
			\measuredangle\prn{\CR\prn{\bbmat \overline{\vek Q}_{j-1}^{\ast}\bbmat \vek 0_{(j-1)L}\\  \\\prod_{i=1}^{j-1}\vek Q_{i}^{(21)} \ebmat\\ \vek 0_{L} \ebmat},\CR\prn{\overline{\vek Q}_{j}^{\ast}\bbmat \vek I_{jL}\\ \vek 0_{L} \ebmat}} &=& \nn\\ \measuredangle\prn{\CR\prn{\bbmat \vek 0_{(j-1)L} \\ \prod_{i=1}^{j-1}\vek Q_{i}^{(21)} \\ \vek 0_{L} \ebmat},\CR\prn{\bbmat \vek I_{L} &&&& \\ &\vek I_{L}&& \\ &&\ddots &\\&&&\vek I_{L}& \\ &&&&\vek Q_{j}^{(11)} \\ \vek 0_{L}  &\vek 0_{L}&\cdots&\vek 0_{L} &\vek Q_{j}^{(12)\ast}\ebmat}}. &&\nn
		\eea
	\normalsize
	We finish the proof by noting that under the assumptions of this paper, we have that $\vek Q_{i}^{(21)}$ is nonsingular
	for $1<i<j-1$ and thus
	\be\nn
		\CR\prn{\bbmat \vek 0_{(j-1)L} \\ \prod_{i=1}^{j-1}\vek Q_{i}^{(21)} \\ \vek 0_{L} \ebmat} = \CR\prn{\bbmat \vek 0_{(j-1)L} \\ \vek I_{L} \\ \vek 0_{L} \ebmat}.
	\ee
	Thus the cosines of the principal angles are the singular values of 
	\footnotesize
	\be\nn
		\bbmat \vek 0_{(j-1)L} & \vek I_{L} & \vek 0_{L} \ebmat\bbmat \vek I_{L} &&&& \\ &\vek I_{L}&& \\ &&\ddots &\\&&&\vek I_{L}& \\ &&&&\vek Q_{j}^{(11)} \\ \vek 0_{L}  &\vek 0_{L}&\cdots&\vek 0_{L} &\vek Q_{j}^{(12)\ast}\ebmat = \bbmat \vek 0_{(j-1)L\times L} & \vek Q_{j}^{(11)} \ebmat
	\ee
	\normalsize
	which are indeed the CS-decomposition cosines, which are the diagonal entries of $\boldsymbol\CC$ from \eqref{eqn.CS-decomp-Q},
	completing the proof.
\eproof
\section{Numerical Examples}\label{section.numerical-examples}
 We constructed two toy 
examples using a matrix considered, e.g., in \cite{B.1991}, 
to demonstrate stagnation properties.  Let
$\vek A_{st}\in\R^{n\times n}$ be defined as the matrix which acts
upon the Euclidean basis as follows,
\be\label{eqn.A-permute}
	\vek A_{st}\vek e_{i}^{[n]} = \bc \vek e_{1}^{[n]}  \mif i=n\\ \vek e_{i+1}^{[n]} \motherwise \ec.
\ee
From this matrix and appropriately chosen right-hand sides, we can generate problems
for which block GMRES is guaranteed to have certain stagnation properties.

In order to obtain some example convergence results in a less non-pathological case,
we also applied block GMRES and FOM to a block diagonal matrix build from $\vek A_{st}$ and 
the {\tt sherman4} matrix from a discretized
oil flow problem, downloaded from the University of Florida Sparse Matrix Library 
\cite{DH.2011}.  The latter matrix is $1104\times 1104$ and nonsymmetric.
\subsection{Total stagnation of block GMRES}
\input{fourStag.tex}

Using the shift matrix $\vek A_{st}$ with $n=200$, we can construct a problem with
four right-hand sides which will stagnate for $50$ iterations before 
converging exactly.  Let the four right-hand sides be the canonical basis vectors
$\vek e_{1}^{[200]}$, $\vek e_{50}^{[200]}$, $\vek e_{100}^{[200]}$, and 
$\vek e_{150}^{[200]}$.  If we let $\vek B\in\R^{200\times 4}$ be the matrix with 
these right-hand sides as columns, we know that
\be\nn
	\vek A^{-1}_{st}\vek B = \bbmat \vek e_{200}^{[200]} & \vek e_{49}^{[200]} & \vek e_{99}^{[200]} & \vek e_{149}^{[200]} \ebmat.
\ee
Due to the stagnating nature of block GMRES for this problem, we compute the
generalized FOM approximation so as to have an iterate at each step.  
The total stagnation for all four right-hand sides can be seen in Figure \ref{figure.fourStag}.

If we arrest the iteration at a stagnating step, e.g., the $40$th step, we can construct
the matrices  $\widetilde{\vek C}_{40}$, $\vek C_{40}$, $\widehat{\vek C}_{40}$, 
$\vek N_{40}$, and $\widehat{\vek N}_{40}$ (all of which are $4\times 4$ matrices) 
to see how such matrices, used to verify theoretical 
results, actually look for a small problem. For the first three matrices, we have the
following,
\be\nn
	\widehat{\vek C}_{40}=\widetilde{\vek C}_{40} = \bbmat 1&&&\\&-1&&\\&&-1&\\&&&-1\ebmat\mand \vek C_{40} = \vek 0_{4},
\ee
and for the last two matrices we have,
\be\nn
	\vek N_{40} = -\vek I_{4} \mand \widehat{\vek N}_{40} = \vek 0_{4}.
\ee
It should be noted that this agrees with what we have proven about block
GMRES in the case of total stagnation in Theorem \ref{thm.H-sing} and
trivially with Theorem \ref{thm.bl-GMRES-FOM-trig-rel}.
\subsection{Partial stagnation/convergence of Block GMRES}
\input{blockStagToy.tex}

In Figure \ref{figure.blockStagToy}, we demonstrate the behavior of block GMRES
and Block FOM applied to a linear system for which block GMRES is guaranteed to 
stagnate but also have earlier convergence for one right-hand side. 
Here, the coefficient matrix is $\vek A_{st}$ defined in \eqref{eqn.A-permute}
for $n=30$.
The block right-hand side $\vek B=\bbmat\vek e_{1}^{[30]}&\vek e_{25}^{[30]}\ebmat$.  
From the definition of $\vek A_{st}$, we have that 
$\vek A_{st}^{-1}\vek B = \bbmat\vek e_{30}^{[30]}&\vek e_{24}^{[30]}\ebmat$.  From this
we see that at iteration $5$, we will achieve exact convergence for the first right-hand
side.  In the absence of replacing the dependent Arnoldi vector with a random one,
the iteration will not produce any improvement for the second right-hand side until
iteration $23$, at which point we again have convergence to the exact solution.
However, in accordance with our block Arnoldi breakdown strategy, we do replace the the dependent basis vector, meaning
we cannot exactly predict stagnation after iteration $5$, though we do see 
near-stagnation until convergence at iteration $15$.

Again, at a particular iterations, we can inspect various
quantities arising which were used in our analysis.  We choose three iterations, $j=5,6,11$, to see what happens at breakdown and
dependent vector replacement.  Indeed we have,
\bea
	\widetilde{\vek C}_{5} &=& \bbmat 1 & 0\\0 &-1\ebmat, \vek C_{5} = \bbmat 0 & 0\\0 &0\ebmat,\mand \widehat{\vek C}_{5}=\bbmat 1 & 0\\0 &-1\ebmat\nn\\
	\widetilde{\vek C}_{6} &=& \bbmat 0 & 0\\-1 &0\ebmat, \vek C_{6} = \bbmat 0 & 0\\-1&0\ebmat,\mand \widehat{\vek C}_{6}=\bbmat 1& 0\\0 &-1\ebmat\nn\\
	\widetilde{\vek C}_{11} &\approx& \bbmat 0.97 & 0\\-0.22 &0\ebmat, \vek C_{11} \approx \bbmat 4.16\times 10^{-17} & 0\\-2.22\times 10^{-16} &0\ebmat,\mand \widehat{\vek C}_{11}=\bbmat 0 & 0\\-1 &0\ebmat\nn
\eea
and we also have 
\bea
	\vek N_{5} &\approx& \bbmat -1.00 & 0\\0 &-1\ebmat \mand \widehat{\vek N}_{5} \approx\bbmat 0 & 0\\0 &0\ebmat\nn\\
	\vek N_{6} &\approx& \bbmat -1.00 & 0\\0 &-1\ebmat \mand \widehat{\vek N}_{6} \approx\bbmat 0 & 1\\0 &0\ebmat\nn\\
	\vek N_{11} &\approx& \bbmat -1.00 & -0.02\\0 &-0.96\ebmat \mand \widehat{\vek N}_{11} \approx\bbmat -0.30 & 0.79\\0 &-1.11\times 10^{-16}\ebmat\nn
\eea
\subsection{A less pathological example with sine computation}
\input{sherman4ToyBlockAnalysis}

To stimulate some slightly more interesting near stagnation behavior, we created a block diagonal matrix in which one 
block is {\tt sherman4} matrix from the University of Florida Sparse Matrix library \cite{DH.2011} and the other block is the shift
matrix $\vek A_{st}$ used in earlier experiments, this time with $n=200$.  
The two right-hand sides are chose to produce perfect stagnation in the shift-matrix block but 
convergence in the {\tt sherman4} block.  Therefore, in the blocks associated to $\vek A_{st}$, the subvectors of 
the right-hand sides were $\vek e_{50}^{[200]}$ and $\vek e_{150}^{[200]}$.  For the {\tt sherman4} matrix, the subvectors
of the right-hand sides were the vector packaged with the matrix and a random vector scaled to have norm
on the order of $10^{7}$.
The exaggerated scaling was done only to produce.
significant convergence prior to stagnation.
In Figure \ref{figure.sherman4ToyBlockAnalysis}, we show the individual $2$-norm
block FOM and block GMRES residual curves as well as the sines from the analysis in Section \ref{section.no-breakdown}.
\section{Conclusions}\label{section.conclusions}
In this paper, we have analyzed the relationship of block GMRES and block
FOM and specifically characterized this relation in the case of block GMRES
stagnation.  These results generalize previous results, particularly those in 
\cite{B.1991} for single vector GMRES and FOM.  We have seen that the relationship
can be a bit more complicated for block methods than in the single-vector method
case due to interaction between approximations for different right-hand sides and 
due to block Arnoldi breakdown.  We close by noting that one
can implement block GMRES so that these sines and cosines are cheaply computable, simply
by following the strategy advocated in \cite{GS.2005} observing that one could implement
a version of block GMRES which also cheaply generates the block FOM approximation.

\section*{Acknowledgment}\label{section.acknowledgment}
We express our gratitude to Mykhaylo Yudytskiy for suggesting the operator
approach for the Moore-Penrose pseudo-inverse. 
We thank Martin Gutknecht for offering many pointers to references on this topic, and Daniel B. Szyld for 
offering helpful comments on the exposition of this paper.
The author also notes that the suggestion to use the CS-decomposition came from Andreas Frommer
while the author was visiting TU-Wuppertal.  

\bibliographystyle{siamplain}
\bibliography{master.bib}
\end{document}

%% file: macros.tex
\def\cal{\mathcal}
\def\norm#1{\left\|#1\right\|}

\def\C{\mathbb{C}}

\def\R{\mathbb{R}}

\def\K{\mathbb{K}}

\def\I{\mathbb{I}}
\def\Cn{\C^n}

\def\BA{{\bf A}}  \def\CA{{\cal A}}
  
  \def\CC{{\cal C}}

  \def\CH{{\cal H}}
  
  \def\CJ{{\cal J}}
  \def\CK{{\cal K}}

  \def\calN{{\cal N}}
  
  \def\CP{{\cal P}}
  \def\CQ{{\cal Q}}
  \def\CR{{\cal R}}
  \def\CS{{\cal S}}
  
  \def\CU{{\cal U}}
  \def\CV{{\cal V}}
  
  \def\CX{{\cal X}}
  \def\CY{{\cal Y}}
  \def\CZ{{\cal Z}}

\def\BAs{\BA{\kern-1.5pt}}

\def\CPs{\CP{\kern-0.8pt}}

\mathcode`\@="8000
{\catcode`\@=\active \gdef@{\mkern1mu}}

\def\mydate{\number\day\ {\ifcase\month \or January\or February\or
              March\or April\or May\or June\or July\or August\or
              September\or October\or November\or December\fi}
\number\year}
\def\vek#1{\mathbf{#1}}

\newcommand{\argmin}[1]{\underset{#1}{\text{{\rm argmin}}}}

\newcommand{\prn}[1]{\left(#1\right)}

\newcommand{\curl}[1]{\left\{#1\right\}}
\newcommand{\ab}[1]{\left|#1\right|}

\newcommand\restr[2]{{
  \left.\kern-\nulldelimiterspace 
  #1 
  \vphantom{\big|} 
  \right|_{#2} 
  }}

\def\Span{\rm span}
\def\rank{\rm rank\ }

\def\be{\begin{equation}}
\def\ee{\end{equation}}
\def\bea{\begin{eqnarray}}
\def\eea{\end{eqnarray}}
\def\nn{\nonumber}
\def\mand{\mbox{\ \ \ and\ \ \ }}
\def\mbut{\mbox{\ \ \ but \ \ }}
\def\mor{\mbox{\ \ \ or\ \ \ }}

\def\mforeach{\mbox{\ \ \ for each\ \ \ }}
\def\mwith{\mbox{\ \ \ with\ \ \ }}
\def\mwhere{\mbox{\ \ \ where\ \ \ }}

\def\mif{\mbox{\ \ \ if\ \ \ }}
\def\motherwise{\mbox{\ \ \ otherwise}}

\def\bbmat{\begin{bmatrix}}
\def\ebmat{\end{bmatrix}}

\def\balg{\begin{algorithm}}
\def\ealg{\end{algorithm}}
\def\bthm{\begin{theorem}}
\def\ethm{\end{theorem}}
\def\blem{\begin{lemma}}
\def\elem{\end{lemma}}
\def\bprop{\begin{proposition}}
\def\eprop{\end{proposition}}
\def\bcor{\begin{corollary}}
\def\ecor{\end{corollary}}
\def\bdefin{\begin{definition}}
\def\edefin{\end{definition}}
\def\bc{\begin{cases}}
\def\ec{\end{cases}}
\newcommand\bproof[1]{\par\addvspace{1ex} \indent\textit{Proof.}\ \ #1}
\def\eproof{\hfill\cvd\linebreak\indent}
\newtheorem{exple}{Example}
\def\bex{\begin{exple}}
\def\eex{\end{exple}}
\newtheorem{assumption}{Assumption}
\def\bass{\begin{assumption}}
\def\eass{\end{assumption}}
\def\diag{{\rm diag}}

\def\cvd{~\vbox{\hrule\hbox{%
  \vrule height1.3ex\hskip0.8ex\vrule}\hrule } }

%% file: fourStag.tex
\begin{figure}[htb]
\hfill
\includegraphics[scale=0.30]{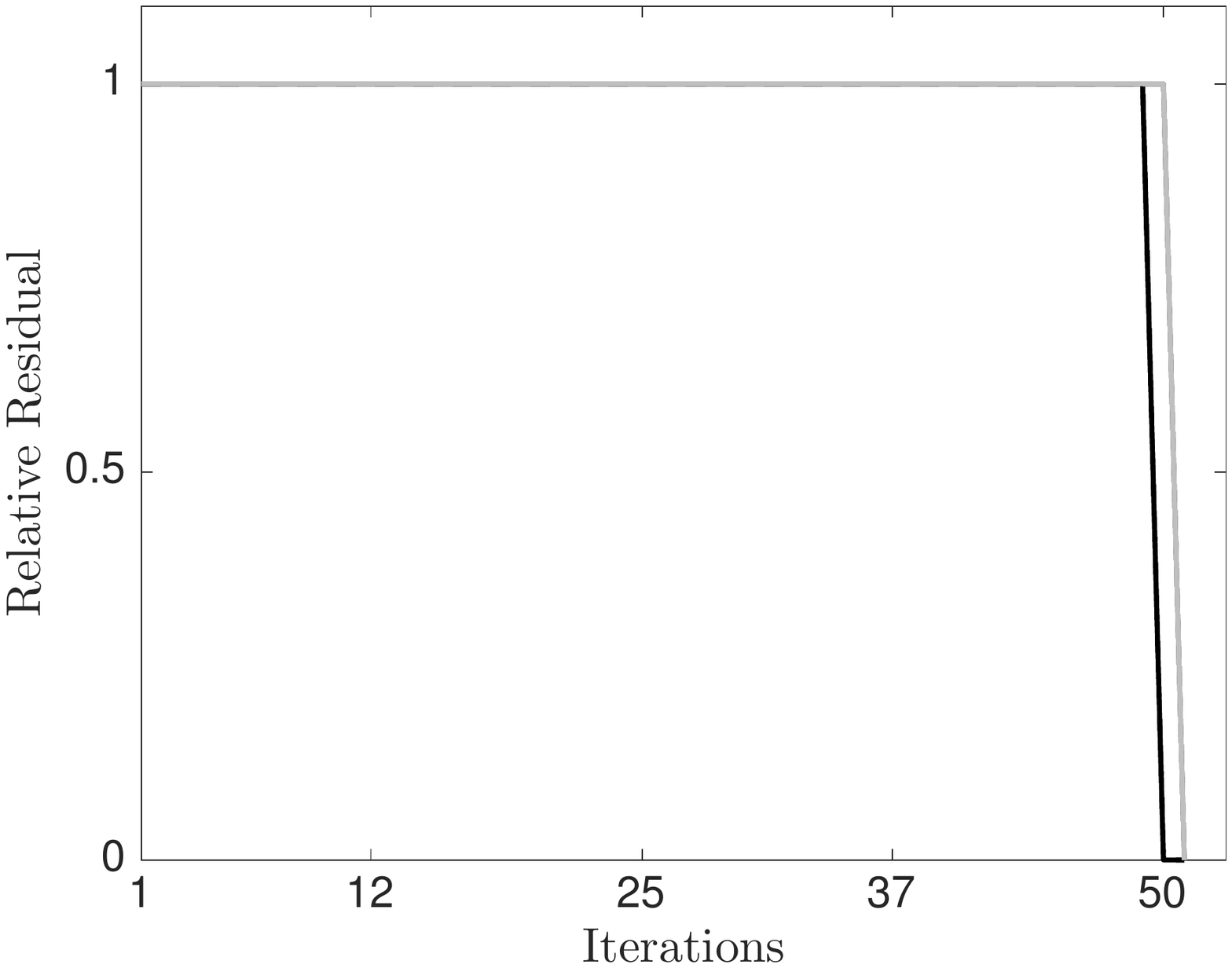}\includegraphics[scale=0.30]{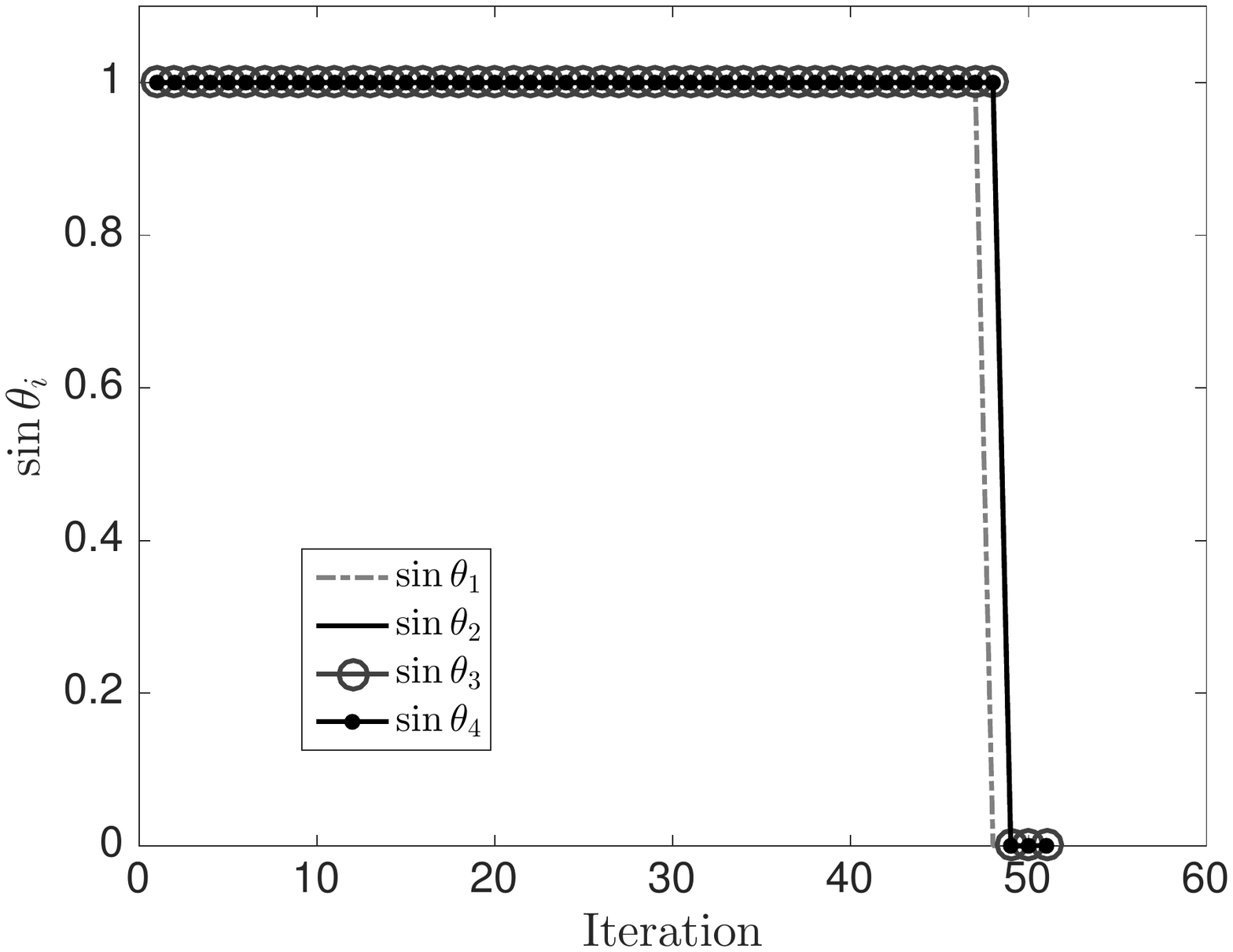}
\hfill
\begin{picture}(0,0)
\end{picture}
\caption{\label{figure.fourStag} Left: Relative two-norm residual curves for stagnating  block GMRES and block FOM for the $200\times 200$ shift matrix for four right-hand sides , namely $\vek e_{1}$, $\vek e_{50}$, $\vek e_{100}$, and $\vek e_{150}$.  The solid and dashed curves  correspond respectively to the block FOM and GMRES residuals, with each shade of gray representing a different  right-hand side.  Similarly, the gray solid and dashed curves, respectively, correspond to the second right-hand side.  Right: Sines of principal angles between $\CR\prn{\vek F_{j-1}^{(G)}}$ and $A\CK_{j}(\vek A, \vek F_{0})$ for each $j$ .}
\end{figure}

%% file: blockStagToy.tex
\begin{figure}[htb]
\hfill
\includegraphics[scale=0.30]{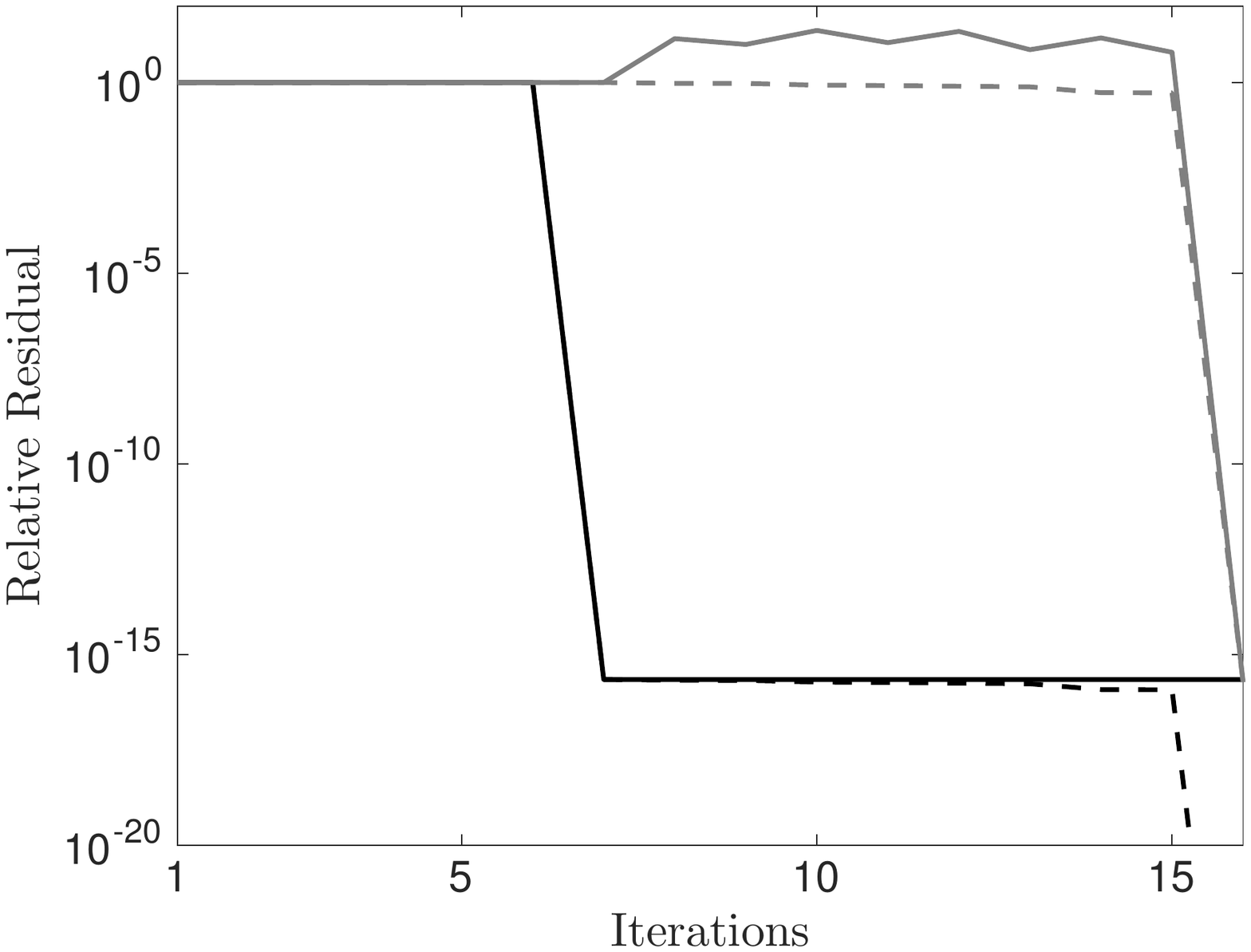}\includegraphics[scale=0.30]{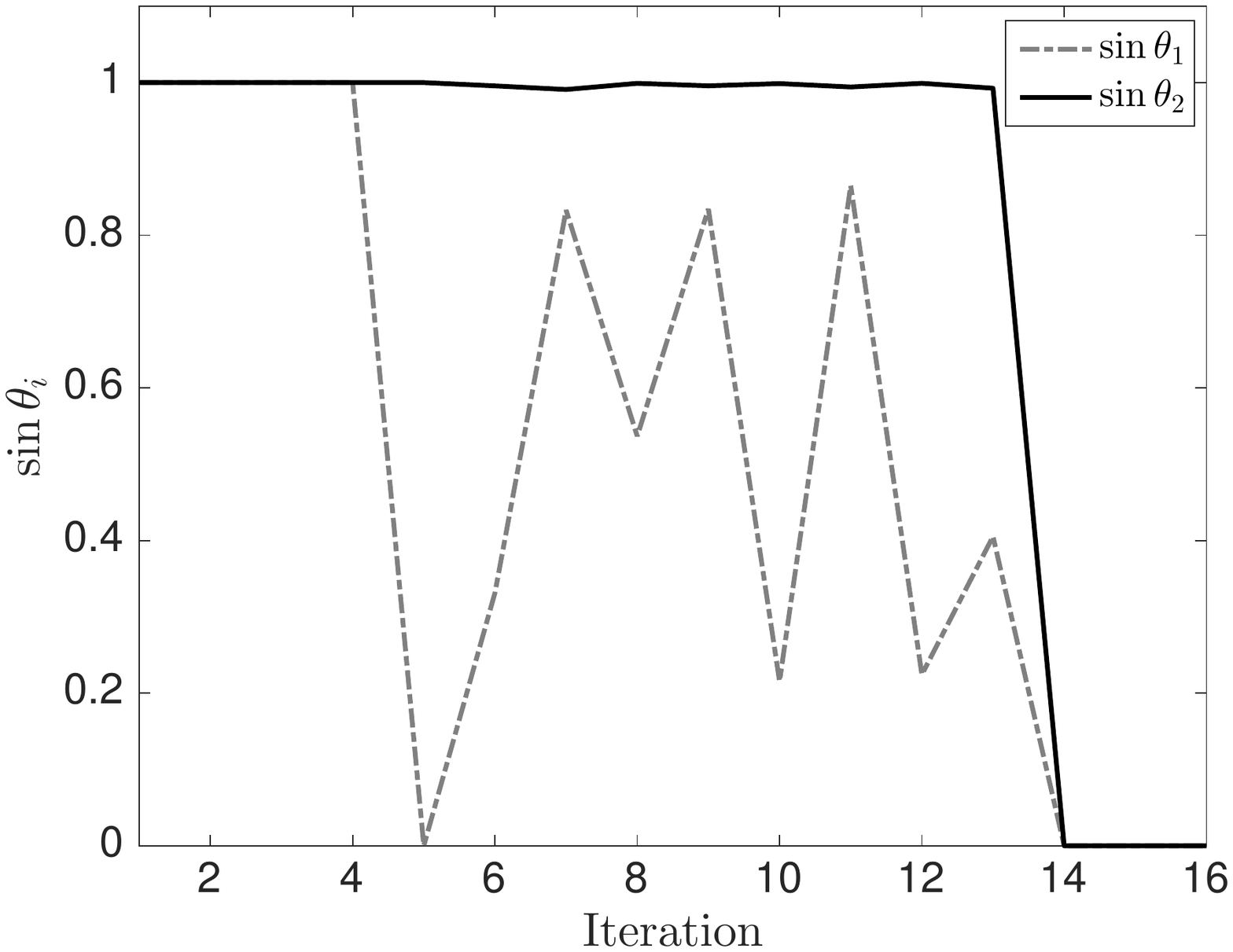}
\hfill
\begin{picture}(0,0)
\end{picture}
\caption{\label{figure.blockStagToy} Left: Relative two-norm residual curves for stagnating  block GMRES and block FOM for the $30\times 30$ shift matrix. The black solid and dashed curves correspond respectively to the  block FOM and GMRES residuals for the first right-hand side.  Similarly,  the gray solid and dashed curves, respectively, correspond to the second right-hand side. Right: Sines of principal angles between $\CR\prn{\vek F_{j-1}^{(G)}}$ and $A\CK_{j}(\vek A, \vek F_{0})$ for each $j$ .}
\end{figure}

%% file: sherman4ToyBlockAnalysis.tex
\begin{figure}[htb]
\hfill
\includegraphics[scale=0.25]{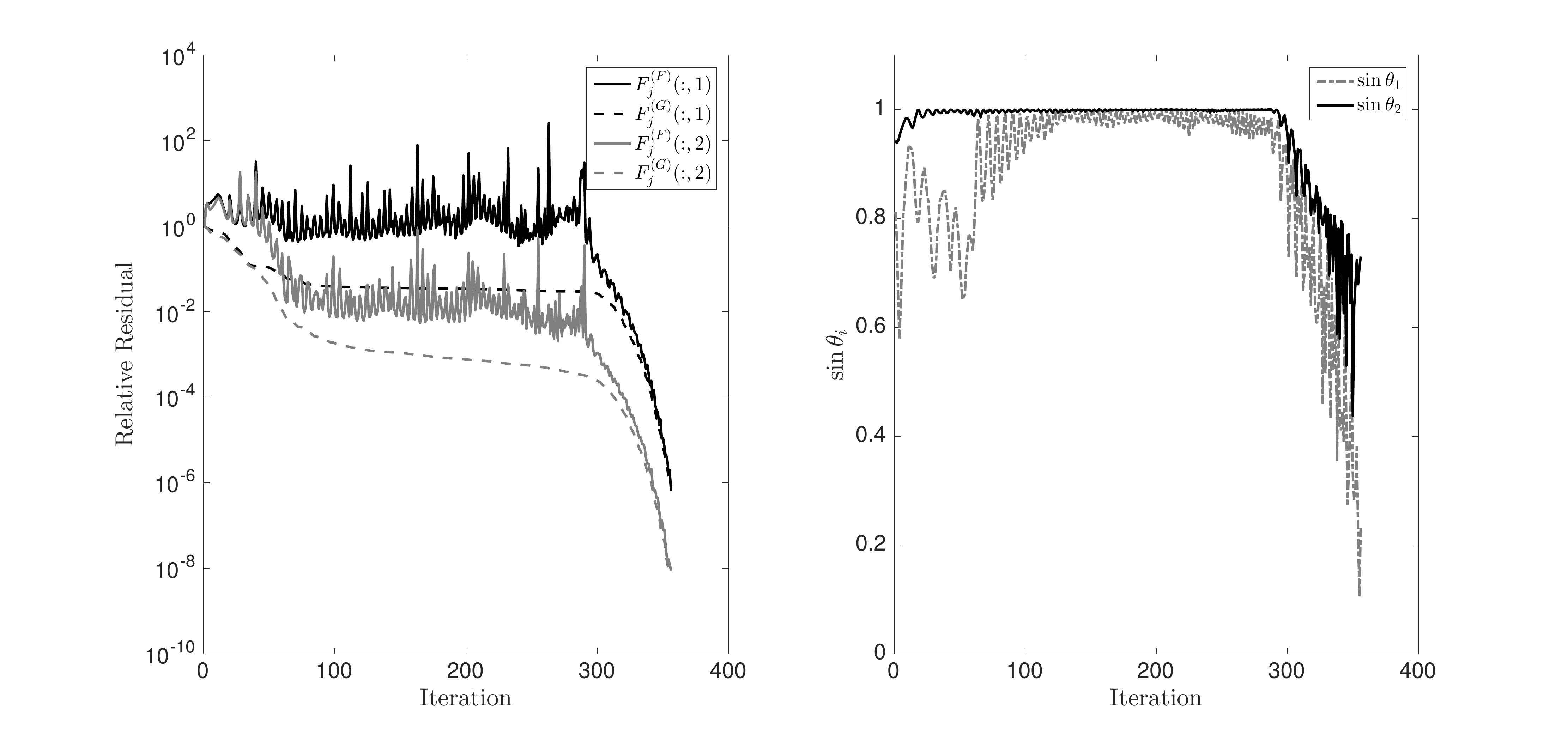}
\hfill
\begin{picture}(0,0)
\end{picture}
\caption{\label{figure.sherman4ToyBlockAnalysis} In the left-hand figure, we have the $2$-norm residual curves of block GMRES and FOM for a linear  system with two right-hand sides using a block diagonal matrix with the {\tt sherman4 } matrix from \cite{DH.2011} as one block  and the shift matrix from the other examples as the other block.  Right-hand sides are chosen to produce wanted near-stagnation. In the right-hand figure, we have the squares of the sines $\curl{s_{1}^{2},s_{2}^{2}}$ coming from the  orthogonal transformations as discussed in the our analysis.}
\end{figure}